\documentclass[12pt]{article}%
\usepackage{amsmath}
\usepackage{amsfonts}
\usepackage{amssymb}
\usepackage{graphicx}%
\newtheorem{theorem}{Theorem}

\newtheorem{proposition}[theorem]{Proposition}

\begin{document}

\author{Diego Dominici \thanks{e-mail: dominicd@newpaltz.edu}\\Department of Mathematics\\State University of New York at New Paltz\\75 S. Manheim Blvd. Suite 9\\New Paltz, NY 12561-2443\\USA\\Phone: (845) 257-2607\\Fax: (845) 257-3571
\and Charles Knessl \thanks{e-mail: knessl@uic.edu}\\Department of Mathematics, Statistics and computer Science\\University of Illinois at Chicago (m/c 249)\\851 South Morgan Street\\Chicago, IL 60607-7045\\USA}
\title{Ray solution of a singularly perturbed elliptic PDE with applications to
communications networks }
\date{}
\maketitle
\begin{abstract}
We analyze a second order, linear, elliptic PDE with mixed boundary
conditions. This problem arose as a limiting case of a Markov-modulated
queueing model for data handling switches in communications networks. We use
singular perturbation methods to analyze the problem. In particular we use the
ray method to solve the PDE in the limit where convection dominates diffusion.
We show that there are both interior and boundary caustics, as well as a cusp
point where two caustics meet, an internal layer, boundary layers and a corner
layer. Our analysis leads to approximate formulas for the queue length (or
buffer content) distribution at the switch.
\end{abstract}


MSC-class: 34E20 (Primary) 60J20 (Secondary)

\section{Introduction}

In a model proposed by Anick, Mitra and Sondhi \cite{AMS}, a buffer receives
messages from $N$ statistically independent and identical information sources,
that asynchronously alternate between exponentially distributed periods in the
``on'' and ``off'' states. While ``on'', a source transmits data at unit rate.
The buffer depletes through an output channel, with a given maximum rate of
transmission $C$. The rate at which a source turns ``on'' is equal to
$\lambda$ and the ``off'' rate is $\mu.$ If $C<N$ the buffer may be non-empty,
and the condition
\[
\frac{\lambda}{\lambda+\mu}N<C
\]
is needed for stability. This simply says that the mean number of ``on''
sources (each transmitting data at unit rate) must be less than the total
transmission capacity of the channel. This model is analyzed exactly in
\cite{AMS}, and the asymptotic limit $N\rightarrow\infty,$ with
\[
\frac{C}{N}=\frac{\lambda}{\lambda+\mu}+O\left(  N^{-\frac{1}{2}}\right)
\]
is studied in \cite{KM}. This limit is referred to as ``heavy traffic''.

Analyzing the steady state joint probability distribution of the number of
active sources and the buffer content, involves solving a system of $N$ linear
ODEs. In heavy traffic this can be simplified to a backward-forward parabolic
PDE of the type in (\ref{parabolic}). This model has the disadvantage of
treating the buffer content as a deterministic fluid.

A modification of this model, which allows for service variability, is as
follows. Again there are $N$ independent and identical sources. When a source
is ``on'' it generates a Poisson arrival stream to a queue. In the ``off''
state no arrivals are generated. The service time distribution is allowed to
be general. The model just described may be called a Markov-modulated M/G/1 queue.

In \cite{KT} it is shown that the joint steady state distribution of the
number of active sources, the queue length and the elapsed service time of the
costumer presently being served satisfies a complicated system of
integro-differential equations. In the heavy traffic limit, where
$N\rightarrow\infty$ and the average arrival rate is close to the mean service
rate, this system may be approximated by the following BVP:
\begin{align}
Df_{yy}+(c-\xi)f_{y}+f_{\xi\xi}+(\xi f)_{\xi}  &  =0,\quad\quad0<y<\infty
,\ -\infty<\xi<\infty\nonumber\\
Df_{y}(0,\xi)+(c-\xi)f(0,\xi)  &  =0,\quad-\infty<\xi<\infty\label{original}\\
\int\limits_{-\infty}^{\infty}\int\limits_{0}^{\infty}f(y,\xi)dyd\xi &
=1.\nonumber
\end{align}
Here the variable $y$ is related to the queue length, $\xi$ corresponds to a
scaled measure of the number of ``on'' sources above their mean value, $c>0 $
is the normalized excess of the service rate over the mean arrival rate, and
$D>0$ measures variability effects in the service time distribution.

The exact solution to (\ref{original}) was analyzed in \cite{KT}. It is not
completely explicit and involves finding one eigenvector of an infinite
matrix, whose elements are complicated expressions involving Laguerre
functions. This (infinite!) eigenvector must be computed numerically. In the
same paper the limit $D\rightarrow\infty$ was considered. Now the matrix
becomes diagonally dominant and much more explicit results can be obtained.

The (highly singular) limit $D\rightarrow0$ was studied in \cite{DK},
resulting in a very complicated asymptotic solution involving contour
integrals of parabolic cylinder and Airy functions. When $D=0$ we see that the
problem (\ref{original}) degenerates into a parabolic one, that is forward
parabolic for $\xi>c$ and backward parabolic for $\xi<c:$%

\begin{align}
(c-\xi)\Im_{y}+\Im_{\xi\xi}+(\xi\Im)_{\xi}  &  =0,\quad\quad0<y<\infty
,\ -\infty<\xi<\infty\nonumber\\
\Im(0,\xi)  &  =0,\quad\ \ \ c<\xi\label{parabolic}\\
\int\limits_{-\infty}^{\infty}\Im(\infty,\xi)d\xi &  =1.\nonumber
\end{align}
Now $\Im$ is a density in $\xi$ and a distribution in $y.$ The problem
(\ref{parabolic}) corresponds to the heavy traffic limit of the fluid model in
\cite{AMS}. Knessl and Morrison \cite{KM} derived the exact solution of
(\ref{parabolic}). The limit $c\rightarrow\infty$ was studied in \cite{K} by
using the saddle point method and in \cite{KK} by using the ray method
\cite{Ke}.

In this paper we will solve (\ref{original}) asymptotically in\ the limit
$c\rightarrow\infty$ by using the ray method, the boundary layer method and
asymptotic matching \cite{BO}. In doing so, we shall analyze no less than
seven different scales, and one more will be briefly discussed in the
conclusion section. The asymptotic structure of (\ref{original}) proves much
more complicated than that of (\ref{parabolic}) in the same limit \cite{KK}.

To analyze (\ref{original}) for large $c$\textit{, }it is convenient to
introduce the new variables $\eta=\xi/c$, $x=y/c$, and the small parameter
$\varepsilon=c^{-2}$. Then (\ref{original}) becomes the following problem for
$F(x,\eta)=\varepsilon^{-1}f(y,\xi)$:
\begin{align}
\varepsilon(DF_{xx}+F_{\eta\eta})+(1-\eta)F_{x}+\eta F_{\eta}+F  &  =0,\quad
x\geq0,\ -\infty<\eta<\infty\nonumber\\
D\varepsilon F_{x}(0,\eta)+(1-\eta)F(0,\eta)  &  =0,\quad-\infty<\eta
<\infty\label{scaled}\\
\int\limits_{-\infty}^{\infty}\int\limits_{0}^{\infty}F(x,\eta)dxd\eta &
=1.\nonumber
\end{align}
The boundary condition together with the normalization condition imply that
the marginal distribution in $\eta$ is the Gaussian
\begin{equation}
\int\limits_{0}^{\infty}F(x,\eta)dx=\frac{1}{\sqrt{2\pi\varepsilon}}%
\exp\left(  -\frac{\eta^{2}}{2\varepsilon}\right)  . \label{marginal}%
\end{equation}
An important quantity to compute is the marginal distribution in the $x$
variable, i.e.,
\begin{equation}
M(x)=\int\limits_{-\infty}^{\infty}F(x,\eta)d\eta. \label{mm}%
\end{equation}

In section 2 we consider the case when $x$ is close to $0$ and $\eta<1;$ this
will be very useful to match with other asymptotic solutions. Section 3 is
dedicated to using the ray method to analyze (\ref{scaled}) for $\varepsilon
\rightarrow0$ with $x,\eta$ fixed. This yields asymptotic solutions in two
main regions separated by the curve $x=\eta-\ln(\eta)-1,\ \eta>1.$ We also
derive boundary layer solutions for $x=O(\varepsilon^{\frac{2}{3}})$ and
$\eta>1,$ $x=O(\varepsilon)$ and $\eta>1$, a corner layer solution in the
neighborhood of the point $(0,1)$ and in section 4 a transition layer solution
along $x=\eta-\ln(\eta)-1.$ We show that all the solutions asymptotically
match to each other in the appropriate limits and also agree with the
approximation found in section 2. In section 5 we summarize and discuss the
main results. In section 6 we check the identity (\ref{marginal}) for
$F(x,\eta)$ and compute the marginal distribution in $x$.

\section{An expansion for small $x$}

To solve (\ref{scaled}) for $\varepsilon$ small, we will first consider the
scaling $x=O(\varepsilon)$. Thus we introduce the variable $v=x/\varepsilon$
and convert (\ref{scaled}) into the problem%

\begin{align}
DF_{vv}+(1-\eta)F_{v}+\varepsilon(\eta F_{\eta}+F)+\varepsilon^{2}F_{\eta
\eta}  &  =0,\quad v\geq0,\ -\infty<\eta<\infty\nonumber\\
DF_{v}(0,\eta)+(1-\eta)F(0,\eta)  &  =0,\quad-\infty<\eta<\infty\label{WKB}\\
\int\limits_{-\infty}^{\infty}\int\limits_{0}^{\infty}F(v,\eta)dvd\eta &
=\frac{1}{\varepsilon}.\nonumber
\end{align}
On this scale (\ref{marginal}) transforms to
\begin{equation}
\int\limits_{0}^{\infty}F(v,\eta)dv=\frac{\varepsilon^{-\frac{3}{2}}}%
{\sqrt{2\pi}}\exp\left(  -\frac{\eta^{2}}{2\varepsilon}\right)  .
\label{WKBmarginal}%
\end{equation}

We consider solutions to (\ref{WKB}) which have the asymptotic form
\begin{equation}
F(v,\eta)\sim\frac{\varepsilon^{-\frac{3}{2}}}{\sqrt{2\pi}}\exp\left(
-\frac{\eta^{2}}{2\varepsilon}\right)  \left[  F^{(0)}(v,\eta)+\sqrt
{\varepsilon}F^{(1)}(v,\eta)+O(\varepsilon)\right]  . \label{WKB1}%
\end{equation}
Substituting (\ref{WKB1}) into (\ref{WKB}) and equating the coefficients of
like powers of $\varepsilon$ we get to leading order the equation
\[
DF_{vv}^{(0)}+(1-\eta)F_{v}^{(0)}=0
\]
with boundary condition
\[
DF_{v}^{(0)}(0,\eta)+(1-\eta)F^{(0)}(0,\eta)=0,\quad-\infty<\eta<\infty.
\]
Solving for $F^{(0)}(v,\eta)$ and taking into account (\ref{WKBmarginal}) we
conclude that%

\[
F^{(0)}(v,\eta)=\frac{1-\eta}{D\sqrt{2\pi}}\exp\left[  -\frac{\eta^{2}%
}{2\varepsilon}-\frac{(1-\eta)v}{D}\right]  ,\quad\eta<1.
\]
We summarize below the main result of this section.

\begin{proposition}
For $x=v\varepsilon=O(\varepsilon)$ the equation (\ref{scaled}) has the
asymptotic solution to leading order
\begin{equation}
F(v,\eta)\sim\frac{1-\eta}{D\sqrt{2\pi}}\exp\left[  -\frac{\eta^{2}%
}{2\varepsilon}-\frac{(1-\eta)v}{D}\right]  ,\quad\eta<1. \label{prop1}%
\end{equation}
\end{proposition}

We see that for $\eta=\overline{\eta}\sqrt{\varepsilon},\ \overline{\eta
}=O(1),$ the solution decouples into a Gaussian in $\overline{\eta}$ times an
exponential function of $v$%
\[
F(v,\overline{\eta})\sim\frac{1}{\sqrt{2\pi}}\exp\left\{  -\frac
{\overline{\eta}^{2}}{2}\right\}  \times\frac{1}{D}\exp\left\{  -\frac{v}%
{D}\right\}  .
\]
Such a decoupling was also observed in \cite{KT}, where (\ref{scaled}) was
analyzed in the limit $D\rightarrow\infty$ with $c$ fixed$.$ The cases where
$x$ is small and $\eta>1$ or $\eta\approx1$ are treated in subsections 3.6 and 3.7$.$

\section{The ray expansion}

Now we consider solutions of (\ref{scaled}) which have the asymptotic form%

\begin{equation}
F(x,\eta)\sim\varepsilon^{\nu_{1}}\exp\left[  \frac{1}{\varepsilon}\Psi
(x,\eta)\right]  K(x,\eta). \label{Ray1}%
\end{equation}

\noindent We substitute (\ref{Ray1}) into (\ref{scaled}) and equate the
coefficients of the lowest power of $\varepsilon$ to get the eikonal equation
for $\Psi:$%

\begin{equation}
D\left(  \Psi_{x}\right)  ^{2}+\left(  \Psi_{\eta}\right)  ^{2}+\eta\left(
\Psi_{\eta}-\Psi_{x}\right)  +\Psi_{x}=0,\quad\Psi_{x}(0,\eta)=\frac{\eta
-1}{D}. \label{eikonal}%
\end{equation}

\noindent Equating the coefficients of the next power of $\varepsilon$ yields
the transport equation for $K:$%

\begin{gather}
DK\Psi_{xx}+K_{x}+2DK_{x}\Psi_{x}+K\Psi_{\eta\eta}+\eta K_{\eta}+2K_{\eta}%
\Psi_{\eta}-\eta K_{x}+K=0,\label{transport}\\
\quad K_{x}(0,\eta)=0.\nonumber
\end{gather}

\subsection{The rays}

We solve (\ref{eikonal}) by introducing the characteristic curves or rays
$\left[  x(t),\ \eta(t)\right]  $, written in terms of a parameter $t$. We
first consider rays starting from the $\eta$-axis, and impose the initial
conditions $\left[  x(0),\eta(0)\right]  =[0,s]$. The characteristic ODEs for
(\ref{eikonal}) are:
\begin{align}
\frac{dx}{dt}  &  =-2D\Psi_{x}+\eta-1,\quad x(0)=0\nonumber\\
\frac{d\eta}{dt}  &  =-2\Psi_{\eta}-\eta,\quad\eta(0)=s\nonumber\\
\frac{d\Psi_{x}}{dt}  &  =0\label{raysys}\\
\frac{d\Psi_{\eta}}{dt}  &  =\Psi_{\eta}-\Psi_{x},\nonumber\\
\frac{d\Psi}{dt}  &  =\Psi_{x}\frac{dx}{dt}+\Psi_{\eta}\frac{d\eta}%
{dt}=-D\left(  \Psi_{x}\right)  ^{2}-\left(  \Psi_{\eta}\right)  ^{2}\nonumber
\end{align}
From (\ref{prop1}) we note that $\Psi(0,\eta)=-\eta^{2}/2$, which implies that
$\Psi\left[  x(0),\eta(0)\right]  =\Psi(0,s)=-s^{2}/2$.

Setting $\Psi_{x}(0,s)=A$, $\Psi_{\eta}(0,s)=B$ and solving (\ref{raysys})
yields:
\begin{align}
x  &  =(A-B)e^{t}-(A+B+s)e^{-t}-(2DA+2A+1)t+2B+s\nonumber\\
\eta &  =(A-B)e^{t}+(A+B+s)e^{-t}-2A\nonumber\\
\Psi_{x}  &  =A\label{ray2}\\
\Psi_{\eta}  &  =(B-A)e^{t}+A\nonumber\\
\Psi &  =-\frac{1}{2}(A-B)^{2}e^{2t}+2A(A-B)e^{t}-A^{2}(D+1)t+AB-\frac{3}%
{2}A^{2}+\frac{1}{2}B^{2}-\frac{s^{2}}{2}.\nonumber
\end{align}
The constants $A,B$ can be determined by evaluating the eikonal equation
(\ref{eikonal}) at $x=0$ (corresponding to $t=0)$, and also using the boundary
condition from (\ref{scaled}). This yields
\[
A=\frac{s-1}{D}\text{ \quad and\quad}B=-s\text{ \ or \ }B=0.
\]

\noindent To decide which value of $B$ is the right one, we take the
derivative of $\Psi$ with respect to $s$ at $t=0$%
\[
-s=\frac{d}{ds}\Psi(0,s)=A\frac{d}{ds}x(0,s)+B\frac{d}{ds}\eta(0,s)=B.
\]

\noindent Replacing $A,B$ in (\ref{ray2}) we get:
\begin{align}
x  &  =e^{t}-1-t-\frac{(D+1)(2t-e^{t})+D+e^{-t}}{D}(s-1)\nonumber\\
\eta &  =e^{t}+\frac{e^{-t}+(D+1)e^{t}-2}{D}(s-1)\label{ray3}\\
\Psi &  =-\frac{1}{2}e^{2t}+\frac{2e^{t}-(D+1)e^{2t}-1}{D}(s-1)\nonumber\\
&  +\frac{-1+[4e^{t}-2(t+1)](D+1)-e^{2t}(D+1)^{2}}{2D^{2}}(s-1)^{2}.\nonumber
\end{align}

\noindent For $t\geq0$ and each value of $s$, the first two equations in
(\ref{ray3}) determine a ray in the $(x,\eta)$ plane, which starts from
$(0,s)$ at $t=0$. For $s=1$ and $\ s=\frac{1}{D+1}$, we can eliminate $t$ from
(\ref{ray3}) and obtain the explicit expressions
\begin{align}
x  &  =X_{0}(\eta)=\eta-\ln(\eta)-1,\quad s=1,\quad\eta\geq
1\label{specialrays}\\
x  &  =\frac{1}{D+1}-\eta-\ln(2-\eta-D\eta),\quad s=\frac{1}{D+1},\quad
\frac{1}{D+1}\leq\eta<\frac{2}{D+1}.\nonumber
\end{align}

\noindent For $s>\frac{1}{D+1}$, we have both $x(t)$ and $\eta(t)$ increasing
for $t>0 $. For $s=\frac{1}{D+1}$, $x(t)$ increases and $\eta(t)$ is
asymptotic to $\frac{2}{D+1}$.

For $s<\frac{1}{D+1}$ the rays ``turn around'' and return to $x=0$ for some
$t^{\ast}>0$, with $x(t^{\ast})=0,\ \eta(t^{\ast})<s$. The maximum value in
$x$ reached bby the ray occurs at $t=t_{x\max}$
\[
t_{x\max}=\ln\left[  \frac{-2sD+D+2-2s+\sqrt{D\left(  4s^{2}D-4sD-8s+4s^{2}%
+D+4\right)  }}{2(1-s-Ds)}\right]
\]%
\begin{figure}
[ptb]
\begin{center}
\rotatebox{270} {\resizebox{11cm}{!}{\includegraphics{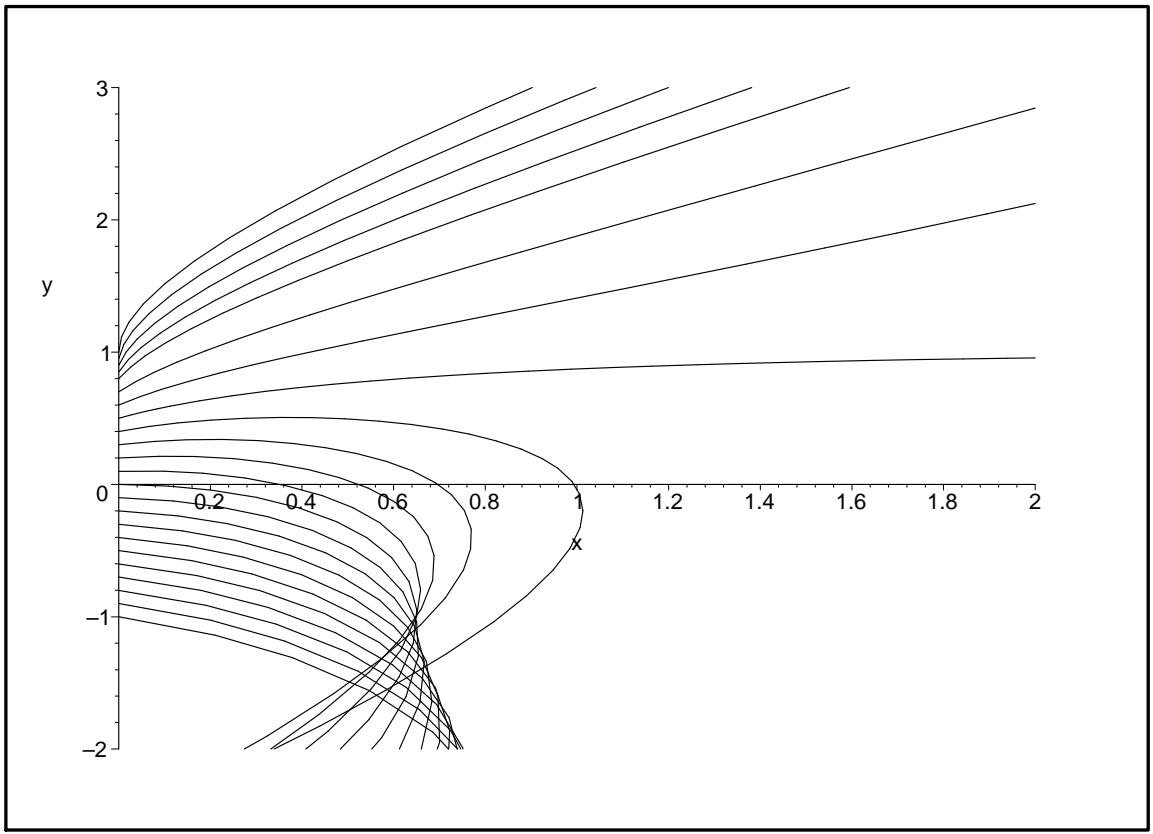}}}
\caption{A sketch of the rays in Region I for $D=1.$}%
\label{Figure1}%
\end{center}
\end{figure}

\noindent For $0<s<\frac{1}{D+1}$ the ray reaches its maximum in $\eta$ at
$t=t_{\eta\max}$
\begin{align*}
t_{\eta\max}  &  =\frac{1}{2}\ln\left[  \frac{1-s}{1-s-Ds}\right]  ,\\
\eta\left(  t_{\eta\max}\right)   &  =2\frac{1-D}{s}+\frac{2s+2D-2-s^{2}%
D-sD^{2}}{s\sqrt{\left(  1-s\right)  \left(  1-s-Ds\right)  }}.
\end{align*}

\noindent For $s\leq0,$ $\eta(t)$ decreases for $0<t<t^{\ast}.$

Solving for $s$ in the $\eta$-equation (\ref{ray3}) yields
\begin{equation}
s=\frac{e^{-t}+e^{t}-2+D\eta}{e^{-t}+(D+1)e^{t}-2} \label{Seta}%
\end{equation}

\noindent and solving in the $x$-equation gives%

\begin{equation}
s=\frac{-e^{t}+e^{-t}+Dt+2t-Dx}{-De^{t}-e^{t}+e^{-t}+2Dt+2t+D}. \label{Sx}%
\end{equation}
Equating (\ref{Seta}) and (\ref{Sx}) we get the implicit equation $R\equiv0$
for the rays, where
\begin{align}
R(x,\eta,t)  &  =\left[  e^{-t}+(D+1)e^{t}-2\right]  x+(3-D\eta-t-Dt-\eta
)e^{t}+(1+t+\eta)e^{-t}\nonumber\\
&  -4-2t+D\eta+2t\eta+2D\eta t. \label{R}%
\end{align}

We sketch several of the rays in Figure \ref{Figure1}. They fill Region I, defined as
\[
\text{Region I }\equiv\text{ }\left\{  x>X_{0}=\eta-\ln(\eta)-1,\quad
\eta>1\right\}  \cup\left\{  x>0,\quad\eta\leq1\right\}  .
\]

\subsection{Caustics and cusps}%

\begin{figure}
[ptb]
\begin{center}
\rotatebox{270} {\resizebox{11cm}{!}{\includegraphics{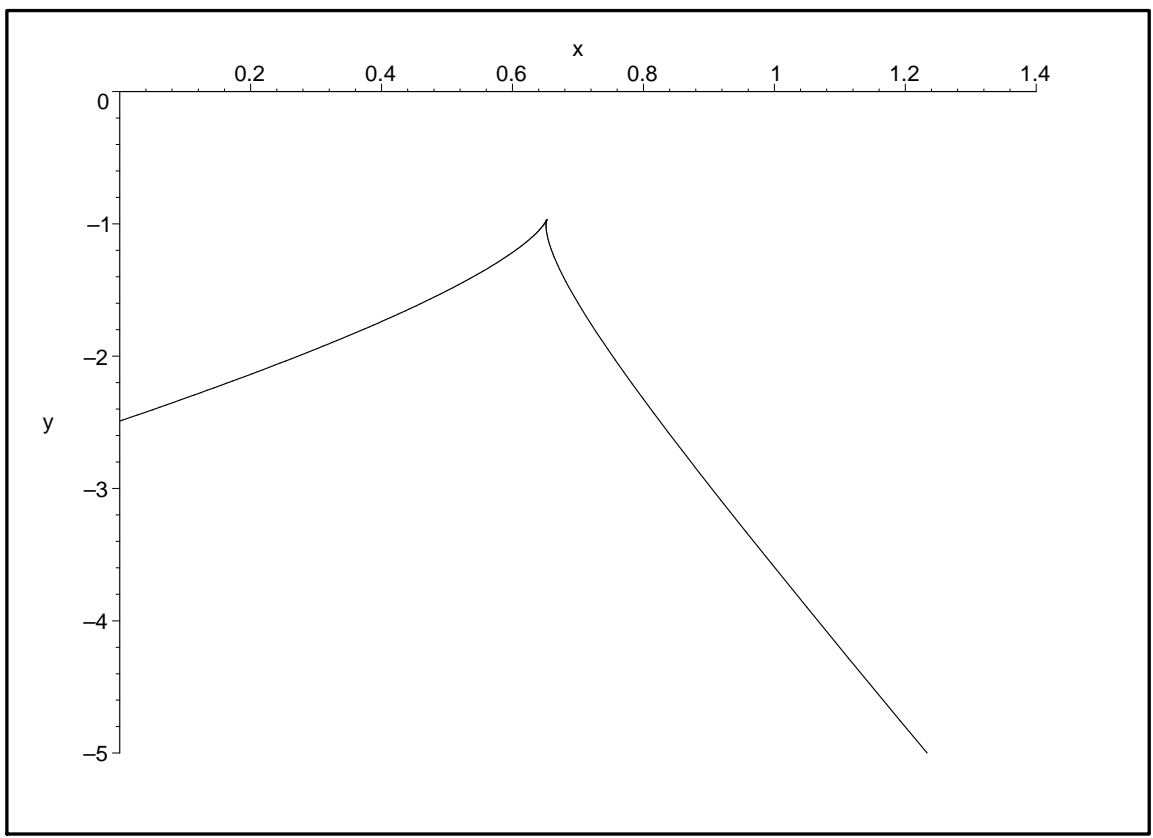}}}
\caption{A sketch of the caustic curves for $D=1.$}%
\label{Figure2}
\end{center}
\end{figure}

The Jacobian of the transformation in (\ref{ray3}) from Cartesian to ray
coordinates is
\begin{align}
J  &  =\frac{dx}{dt}\frac{d\eta}{ds}-\frac{dx}{ds}\frac{d\eta}{dt}\label{J}\\
&  =\left[  2(t-2)(s-1)D^{-2}+(-2t-5s+4ts+2)D^{-1}-s+2ts+1\right]
e^{t}\nonumber\\
&  +\left[  -2(t+2)(s-1)D^{-2}+(2t-2ts+2-3s)D^{-1}\right]  e^{-t}\nonumber\\
&  +8(s-1)D^{-2}+4(2s-1)D^{-1}\nonumber
\end{align}
When $J=0$ we can solve for $s$ as a function of $t,\ S_{0}=s\left|
_{J=0}\right.  $%
\begin{equation}
S_{0}=\frac{(-2D-D^{2}-4+2Dt+2t)e^{2t}+4(D+2)e^{t}-2(2+D+Dt+t)}{(-D^{2}%
-5D-4+2t+4Dt+2tD^{2})e^{2t}+8(D+1)e^{t}-3D-4-2t-2Dt}. \label{SJ0}%
\end{equation}

The equation for the \textit{caustic(s),} i.e., the points in the $(x,\eta)$
plane at which the Jacobian is zero, can be given in parametric form. We
replace $s$ by $S_{0}$ in the equation of the rays, and let $x_{ca}%
=x(t,S_{0}),\ \eta_{ca}=\eta(t,S_{0})$:
\begin{align}
x_{ca}  &  =\left[  -(D+1)^{2}e^{3t}+(2D^{2}t^{2}-3tD+D^{2}t+2t^{2}%
-4t+D^{2}+4t^{2}D+6D+8)e^{2t}\right. \nonumber\\
&  \left.  -2(3D+7)e^{t}-e^{-t}+2(D+1)t^{2}+(3D+4)t+2(D+4)\right] \label{Xc}\\
&  /\left[  (2D^{2}t+4Dt-4+2t-D^{2}-5D)e^{2t}+8(D+1)e^{t}%
-(3D+4)-2(D+1)t\right] \nonumber
\end{align}%
\begin{align}
\eta_{ca}  &  =\left[  -(D+1)^{2}e^{3t}+2(2tD+2t+2D-1)e^{2t}%
+2(4-2t-2tD-D)e^{t}+e^{-t}-6\right] \nonumber\\
&  /\left[  (2D^{2}t+4Dt-4+2t-D^{2}-5D)e^{2t}+8(D+1)e^{t}%
-(3D+4)-2(D+1)t\right]  \label{etaC}%
\end{align}

In Figure \ref{Figure2} we sketch the caustic curves for $D=1$. There is also a cusp
where the two caustics meet. Our numerical studies show that the basic
structure (i.e., the two caustics coming together as a cusp) occurs for all
$D>0.$ \ 

\begin{figure}
[ptb]
\begin{center}
\rotatebox{270} {\resizebox{11cm}{!}{\includegraphics{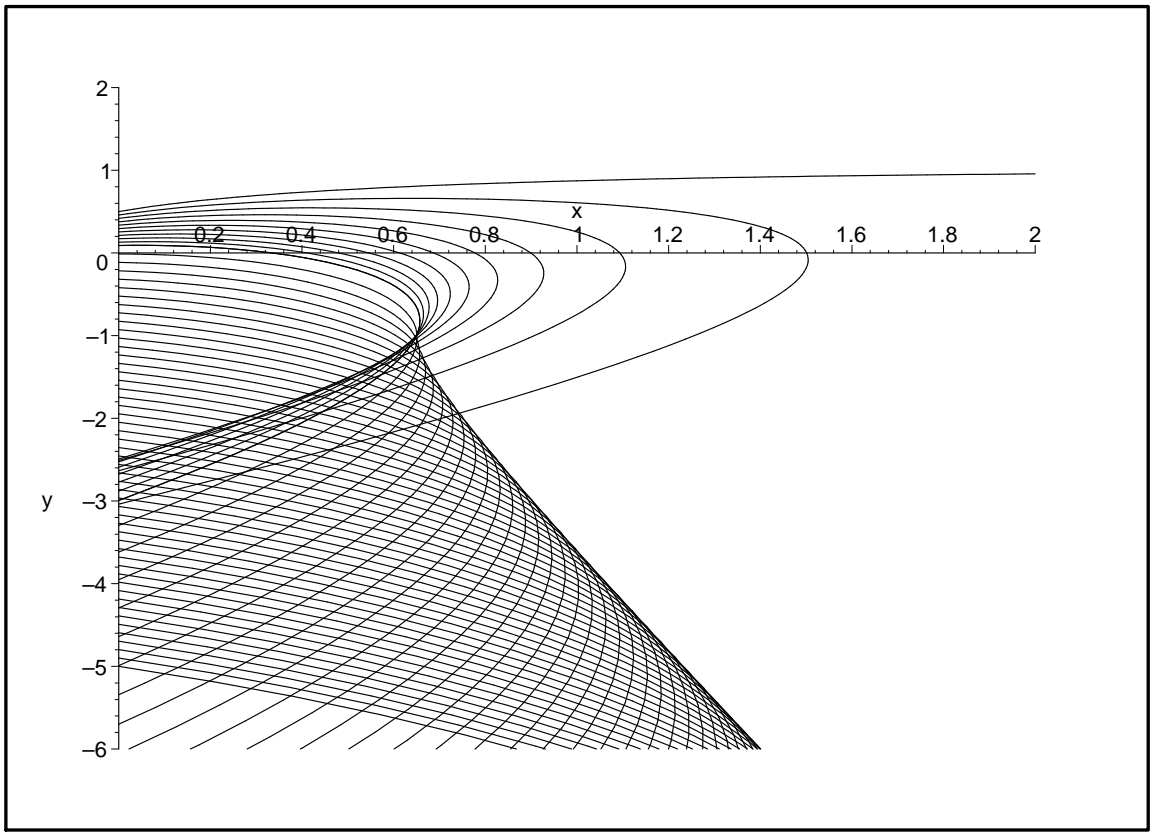}}}
\caption{A sketch of the rays in Region I for $D=1.$}%
\label{Figure3}%
\end{center}
\end{figure}

Outside the caustic region, the correspondence between $(t,s)$ and $(x,\eta)$
is one-to-one. When we are exactly on the caustic curves, the correspondence
is two-to-one, and inside the region bounded by the two caustics it is
three-to-one. In Figure \ref{Figure3} we sketch more densely the rays for $D=1$ to
indicate this correspondence. The evaluation of (\ref{Ray1}) near caustics and
cusps is discussed in more detail in section 5.%

\subsection{The transport equation}

Now we shall solve the transport equation (\ref{transport}) by using
(\ref{raysys}) to write it as an ODE along a ray:
\begin{equation}
\frac{dK}{dt}=(D\Psi_{xx}+\Psi_{\eta\eta}+1)K. \label{transport1}%
\end{equation}
After some algebra, we can show that
\[
D\Psi_{xx}+\Psi_{\eta\eta}+1=\frac{1}{2}-\frac{1}{2J}\frac{dJ}{dt}%
\]
and hence
\[
K(x,\eta)=k(s)\frac{e^{\frac{t}{2}}}{\sqrt{J}}.
\]

To determine $k(s)$ we evaluate the previous result at $t=0:$%
\[
K(0,s)=k(s)\frac{1}{\sqrt{1-s}}.
\]
Using the approximation (\ref{prop1}) and the fact that $s=\eta$ at $t=0$, we
get
\[
k(s)=\frac{1}{\sqrt{2\pi}D}(1-s)^{\frac{3}{2}},\quad s<1\quad\text{ and \ }%
\nu_{1}=-\frac{3}{2}.
\]

The same result can be obtained by using the BC $K_{x}(0,\eta)=0$ in
(\ref{transport}) and fixing the multiplicative constant by normalization. So
far we have determined $\Psi$ and $K$ only for $s<1$. Thus we divide the
half-plane $x\geq0,\ -\infty<\eta<\infty$ into two parts. The portion filled
by the rays for $s<1$ we call Region I and the remainder of the half-plane we
call Region II. The latter is a shadow of the rays (see also Figure 3.1).

To summarize, we have established the following.

\begin{proposition}
The solution of (\ref{scaled}) in Region I is asymptotically given by
\[
F(x,\eta)\sim\varepsilon^{-\frac{3}{2}}K(x,\eta)\exp\left[  \frac
{1}{\varepsilon}\Psi(x,\eta)\right]
\]
where
\begin{align}
K(x,\eta)  &  =\frac{1}{\sqrt{2\pi}}(1-s)^{\frac{3}{2}}\frac{e^{\frac{t}{2}}%
}{\sqrt{J(t,s)}}\nonumber\\
\Psi(x,\eta)  &  =-\frac{1}{2}e^{2t}+\frac{2e^{t}-(D+1)e^{2t}-1}%
{D}(s-1)\label{prop2}\\
&  +\frac{-1+[4e^{t}-2(t+1)](D+1)-e^{2t}(D+1)^{2}}{2D^{2}}(s-1)^{2},\nonumber
\end{align}
$(x,\eta)$ is related to $(t,s)$ by (\ref{ray3}) and $J(t,s)$ is defined by
(\ref{J}).
\end{proposition}

\subsection{Region II}

For this region, we consider solutions of (\ref{scaled}) which have the
asymptotic form
\[
F(x,\eta)\sim\varepsilon^{\nu_{2}}\exp\left[  \frac{1}{\varepsilon}\Phi
(x,\eta)+\frac{1}{\varepsilon^{\frac{1}{3}}}\Gamma(x,\eta)\right]  L(x,\eta)
\]
The term $\varepsilon^{-\frac{1}{3}}\Gamma(x,\eta)$ in the exponent must be
included in order for the expansion to asymptotically match those valid for
small $x$ and $\eta>1$, which we construct later.

It follows that $\Phi$ satisfies (\ref{eikonal}), $L$ satisfies
(\ref{transport}) and for $\Gamma$ we get the following PDE
\begin{equation}
(\eta-1-2D\Phi_{x})\Gamma_{x}-(2\Phi_{\eta}+\eta)\Gamma_{\eta}=0,
\label{Gamma}%
\end{equation}
which is equivalent to $\frac{d\Gamma}{d\tau}=0.$ Thus we conclude that
$\Gamma$ is a function of $\sigma$ only and write $\Gamma(x,\eta
)=\Gamma(\sigma)$. Here $(\tau,\sigma)$ are the new parameters for the ray
which apply in Region II. Thus a ray starts at $\tau=0$ from $\eta=\sigma>1$
and enters the domain for $\tau>0.$

The solutions of the characteristic equations are:
\begin{align}
x  &  =(b-a)e^{\tau}+(a+b-\sigma)e^{-\tau}+[2a(D+1)-1]\tau-2b+\sigma
\nonumber\\
\eta &  =(b-a)e^{\tau}-(a+b-\sigma)e^{-\tau}+2a\nonumber\\
\Phi_{x}  &  =-a\label{rayII1}\\
\Phi_{\eta}  &  =(a-b)e^{\tau}-a\nonumber\\
\Phi &  =-a^{2}(D+1)\tau+2a\left(  a-b\right)  \left(  e^{\tau}-1\right)
-\frac{1}{2}\left(  a-b\right)  ^{2}\left(  e^{2\tau}-1\right)  +\Phi
_{0}(\sigma).\nonumber
\end{align}
Here $\Phi_{0}(\sigma)$ is the value of $\Phi$ at $\tau=0$, which corresponds
to the $\eta$-axis for $\eta>1.$

Since from the result for Region I $\frac{dx}{dt}=\frac{(s-1)}{D}=0\ $for
$s=1,$ we impose the condition $\frac{dx}{d\tau}(0,\sigma)=0$ for all
$\sigma>1.$ This means that the boundary $x=0$ will be a caustic curve for
$\eta>1.$ Then $a$ has the value
\begin{equation}
a(\sigma)=\frac{1-\sigma}{2D}. \label{a}%
\end{equation}
Evaluating (\ref{eikonal}) at $x=0$ we get
\begin{equation}
Da^{2}+b^{2}+\sigma(b-a)+a. \label{strip}%
\end{equation}
Using (\ref{a}) in (\ref{strip}) and solving for $b$ we find that%

\begin{equation}
b=\frac{\sigma}{2}\pm\frac{\sqrt{\beta(\sigma)}}{2\sqrt{D}},\quad\beta
(\sigma)=D\sigma^{2}+\left(  \sigma-1\right)  ^{2}. \label{beta}%
\end{equation}

For small $\tau$ we get from (\ref{rayII1}) and (\ref{a})
\[
x\sim\left(  b-\frac{\sigma}{2}\right)  \tau^{2},\quad\tau\rightarrow0
\]
and this implies that the solution $b=\frac{\sigma}{2}-\frac{\sqrt
{\beta(\sigma)}}{2\sqrt{D}}$ must be rejected, in order that the rays enter
the domain $x\geq0,$ as $\tau$ increases. Hence,
\begin{equation}
b(\sigma)=\frac{\sigma}{2}+\frac{\sqrt{\beta(\sigma)}}{2\sqrt{D}}. \label{b}%
\end{equation}

To find $\Phi_{0}(\sigma)$ we impose the continuity condition $\Phi
_{0}(1)=\Psi(0,1)=-\frac{1}{2}.$ Since
\[
\frac{d}{d\sigma}\Phi(0,\sigma)=-a\frac{d}{d\sigma}x(0,\sigma)-b\frac
{d}{d\sigma}\eta(0,\sigma)=-b
\]
we conclude that
\begin{align}
\Phi_{0}(\sigma)  &  =-\frac{1}{2}-\int\limits_{1}^{\sigma}b\left(  u\right)
du\nonumber\\
&  =-\frac{1}{4}-\frac{\sigma^{2}}{4}-\frac{1}{4\sqrt{D}}\left\{  \left[
\sigma-\frac{1}{D+1}\right]  \sqrt{\beta(\sigma)}\right. \label{phi00}\\
&  +\frac{D}{(D+1)^{\frac{3}{2}}}\operatorname{arcsinh}\left[  \frac
{(D+1)\sigma-1}{\sqrt{D}}\right] \nonumber\\
&  \left.  -\frac{D^{\frac{3}{2}}\sigma}{(D+1)}-\frac{D}{(D+1)^{\frac{3}{2}}%
}\operatorname{arcsinh}\left[  \sqrt{D}\right]  \right\}  .\nonumber
\end{align}

As before, the transport equation (\ref{transport}) can be solved to obtain
\begin{equation}
L(\tau,\sigma)=L_{0}(\sigma)\frac{e^{\frac{\tau}{2}}}{\sqrt{\widetilde{J}}}
\label{L}%
\end{equation}
where%

\begin{align}
\widetilde{J}  &  =\frac{dx}{d\tau}\frac{d\eta}{d\sigma}-\frac{dx}{d\sigma
}\frac{d\eta}{d\tau}\nonumber\\
&  =\left\{  \left[  -\sigma+1+\frac{1}{2}\tau(\sigma-1)\right]
D^{-2}+\left[  \frac{1}{2}\sqrt{\beta(\sigma)}(\tau-1)\right]  D^{-\frac{3}%
{2}}\right. \nonumber\\
&  \left.  +\left(  -\sigma-\frac{1}{2}\tau+\tau\sigma\right)  D^{-1}+\frac
{1}{2}\tau\sqrt{\beta(\sigma)}D^{-\frac{1}{2}}+\frac{1}{2}\tau\sigma\right\}
e^{\tau}\label{J1}\\
&  +\left\{  \left(  \frac{1}{2}\tau+1\right)  (1-\sigma)D^{-2}+\left[
\frac{1}{2}\sqrt{\beta(\sigma)}(\tau+1)\right]  D^{-\frac{3}{2}}\right.
\nonumber\\
&  \left.  +\left(  -\sigma+\frac{1}{2}\tau-\tau\sigma\right)  D^{-1}+\frac
{1}{2}\tau\sqrt{\beta(\sigma)}D^{-\frac{1}{2}}-\frac{1}{2}\tau\sigma\right\}
e^{-\tau}\nonumber\\
&  +2(\sigma-1)D^{-2}+2\sigma D^{-1}.\nonumber
\end{align}

\noindent In Region II $\widetilde{J}=0$ only for $\tau=0.$ To determine
$L_{0}(\sigma)$ and $\Gamma(\sigma)$ we shall analyze the problem for small $x
$, and we will find that not one, but two boundary layer expansions are needed
to satisfy the boundary conditions (\ref{scaled}) in this region.

\subsection{Approximation for $x=O(\varepsilon^{\frac{2}{3}}),\ \eta>1$ (inner solution)}

We introduce the stretched variable $\mu=\varepsilon^{-\frac{2}{3}}x$, and
transform (\ref{scaled}) into
\begin{equation}
(1-\eta)F_{\mu}+\varepsilon^{\frac{1}{3}}DF_{\mu\mu}+\varepsilon^{\frac{2}{3}%
}(\eta F_{\eta}+F)+\varepsilon^{\frac{5}{3}}F_{\eta\eta}=0. \label{scaled2/3}%
\end{equation}
We represent $F$ in the asymptotic form
\begin{equation}
F\sim\varepsilon^{\nu_{3}}\exp\left\{  \varepsilon^{-1}\Phi_{0}(\eta
)+\varepsilon^{-\frac{1}{3}}\left[  \frac{\eta-1}{2D}\mu+\Gamma\left(
\eta\right)  \right]  \right\}  \left[  R_{0}(\mu,\eta)+\varepsilon^{\frac
{1}{3}}R_{1}(\mu,\eta)\right]  \label{F2/3}%
\end{equation}
which when inserted into (\ref{scaled2/3}) give the following PDEs for $R_{0}
$, $R_{1}:$%
\begin{equation}
2D^{2}\frac{\partial^{2}R_{0}}{\partial\mu^{2}}+\left[  2\Phi_{0}^{\prime
}(\eta)+\eta\right]  \left[  2D\Gamma^{\prime}(\eta)+\mu\right]  R_{0}=0
\label{R0}%
\end{equation}%

\begin{align}
0  &  =2D^{2}\frac{\partial^{2}R_{1}}{\partial\mu^{2}}+\left[  2\Phi
_{0}^{\prime}(\eta)+\eta\right]  \left[  2D\Gamma^{\prime}(\eta)+\mu\right]
R_{1}\label{R1}\\
&  +2D\left\{  \left[  2\Phi_{0}^{\prime}(\eta)+\eta\right]  \frac{\partial
R_{0}}{\partial\eta}+\left[  \Phi_{0}^{\prime\prime}(\eta)+1\right]
R_{0}\right\} \nonumber
\end{align}

Solving (\ref{R0}) we get
\begin{equation}
R_{0}=C_{1}(\eta)\operatorname{Ai}\left\{  2^{-\frac{1}{3}}D^{-\frac{5}{6}%
}\beta\left(  \eta\right)  ^{\frac{1}{6}}\left[  \mu+2D\Gamma^{\prime}%
(\eta)\right]  \right\}  \label{R0sol}%
\end{equation}
where $\operatorname{Ai}(\cdot)$ denotes the Airy function and $\beta\left(
\eta\right)  $ is given by (\ref{beta}). Using (\ref{R0sol}) into (\ref{R1})
and solving for $R_{1}$ we obtain
\begin{align}
R_{1}  &  =\frac{1}{24}2^{\frac{2}{3}}D^{\frac{1}{6}}\beta(\eta)^{-\frac{5}%
{6}}C_{1}(\eta)\beta^{\prime}(\eta)\overline{\mu}^{2}\operatorname{Ai}%
(\overline{\mu})+\left[  2D\beta(\eta)\right]  ^{\frac{1}{3}}C_{1}(\eta
)\Gamma^{\prime\prime}(\eta)\overline{\mu}\operatorname{Ai}(\overline{\mu
})+\nonumber\\
&  +\left\{  2^{\frac{2}{3}}\left[  D\beta(\eta)\right]  ^{\frac{1}{6}}%
C_{1}^{\prime}(\eta)-2^{-\frac{1}{3}}D^{\frac{2}{3}}\beta(\eta)^{-\frac{1}{3}%
}C_{1}(\eta)+\frac{1}{3}2^{\frac{2}{3}}D^{\frac{1}{6}}\beta(\eta)^{-\frac
{5}{6}}\alpha(\eta)C_{1}(\eta)\right\}  \operatorname{Ai}^{\prime}%
(\overline{\mu})\label{R1sol}\\
&  +C_{2}(\eta)\operatorname{Ai}(\overline{\mu})\nonumber
\end{align}
with
\[
\overline{\mu}=2^{-\frac{1}{3}}D^{-\frac{5}{6}}\beta\left(  \eta\right)
^{\frac{1}{6}}\left[  \mu+2D\Gamma^{\prime}(\eta)\right]  ,\quad
\]%
\begin{equation}
\alpha(\eta)=\left(  D+1\right)  \eta-1. \label{alpha}%
\end{equation}
The function $C_{1}(\eta)$ will be determined below. This solution can't
satisfy the boundary condition (\ref{scaled}), and thus we require another
boundary layer expansion, where $x=o(\varepsilon^{\frac{2}{3}})$.

\subsection{Approximation for $x=O(\varepsilon),\ \eta>1$ (inner-inner solution)}

We introduce the variable $v=x/\varepsilon,$ and transform (\ref{scaled}) to
\begin{align}
DF_{vv}+(1-\eta)F_{v}+\varepsilon(\eta F_{\eta}+F)+\varepsilon^{2}F_{\eta
\eta}  &  =0\label{scaled1}\\
DF_{v}(0,\eta)+(1-\eta)F(0,\eta)  &  =0.\nonumber
\end{align}

We seek solutions of the form
\begin{equation}
F\sim\varepsilon^{\nu_{4}}\exp\left\{  \frac{1}{\varepsilon}\Phi_{0}%
(\eta)+\frac{1}{\varepsilon^{\frac{1}{3}}}\Gamma(\eta)+\frac{1}{2}\frac
{\eta-1}{D}v\right\}  W(v,\eta). \label{F1}%
\end{equation}
Using (\ref{F1}) in (\ref{scaled1}) and taking into account that
\[
\Phi_{0}^{\prime}(\eta)=-b(\eta)=-\left(  \frac{\eta}{2}+\frac{\sqrt{D\eta
^{2}+\left(  \eta-1\right)  ^{2}}}{2\sqrt{D}}\right)
\]
yields
\begin{align*}
DW_{vv}+(1-\eta)W_{v}+\Phi_{0}^{\prime}(\eta)\left[  \Phi_{0}^{\prime}%
(\eta)+\eta\right]  W  &  =0\\
2DW_{v}(0,\eta)+(1-\eta)W(0,\eta)  &  =0
\end{align*}
whose general solution is
\begin{equation}
W(v,\eta)=w(\eta)\left[  \frac{1}{2D}(\eta-1)v+1\right]  \label{Wsol}%
\end{equation}

The next step will be finding a corner layer solution valid in a neighborhood
of the point $(0,1),$ that matches to both the approximation (\ref{prop1}) and
the inner-inner solution. This will allow us to determine $\nu_{4}$ and
$w(\eta)$ explicitly.

\subsection{Corner layer}

Let us first write $F(x,\eta)=\varepsilon^{\nu_{5}}\exp\left(  -\frac{\eta
^{2}}{2\varepsilon}\right)  \overline{G}(x,\eta)$, which transforms
(\ref{scaled}) into
\begin{align}
D\varepsilon\overline{G}_{xx}-\eta\overline{G}_{\eta}+\varepsilon\overline
{G}_{\eta\eta}+(1-\eta)\overline{G}_{x}  &  =0\label{scaledcorner}\\
D\varepsilon\overline{G}_{x}(0,\eta)+(1-\eta)\overline{G}(0,\eta)  &
=0.\nonumber
\end{align}
Then we introduce the stretched variables $\mu=\varepsilon^{-\frac{2}{3}}x$
and $\gamma=\varepsilon^{-\frac{1}{3}}\left(  \eta-1\right)  ,$ and
(\ref{scaledcorner}) becomes
\begin{align}
\varepsilon^{\frac{2}{3}}\overline{G}_{\gamma\gamma}-\varepsilon^{\frac{1}{3}%
}\gamma\overline{G}_{\gamma}+D\overline{G}_{\mu\mu}-\gamma\overline{G}_{\mu
}-\overline{G}_{\gamma}  &  =0\label{scaledcorner1}\\
D\overline{G}_{\mu}(0,\gamma)-\gamma\overline{G}(0,\gamma)  &  =0.\nonumber
\end{align}
To leading order $\overline{G}(\mu,\gamma)\sim G(\mu,\gamma)$ where
\begin{align}
DG_{\mu\mu}-\gamma G_{\mu}-G_{\gamma}  &  =0\label{scaledcornerQ}\\
DG_{\mu}(0,\gamma)-\gamma G(0,\gamma)  &  =0.\nonumber
\end{align}

The solution to (\ref{scaledcornerQ}) matches to (\ref{prop1}) (with $\mu=0) $
if
\begin{equation}
\varepsilon^{\nu_{5}}G(0,\gamma)\sim\frac{1-\eta}{\sqrt{2\pi}D}\varepsilon
^{-\frac{3}{2}}=-\frac{\gamma}{\sqrt{2\pi}D}\varepsilon^{-\frac{7}{6}}%
,\quad\gamma\rightarrow-\infty\label{Glimit}%
\end{equation}
so that $\nu_{5}=-\frac{7}{6}.$ In \cite{knesslcorner} an explicitly solution
to (\ref{scaledcornerQ}) and (\ref{Glimit}) was obtained \footnote{The
function $Q(X,T)$ in \cite{knesslcorner} is related to $G(\mu,\gamma)$ by
\[
G(\mu,\gamma)=\frac{1}{\sqrt{2\pi}}(2D)^{-\frac{2}{3}}Q\left[  (2D)^{-\frac
{2}{3}}\mu,(2D)^{-\frac{1}{3}}\gamma\right]  ,\quad r_{0}=2^{\frac{1}{3}}%
\beta_{\ast}.
\]
}, with
\[
G(\mu,\gamma)=\frac{\exp\left\{  \frac{1}{\varepsilon}\left[  \frac{\mu\gamma
}{2D}-\frac{\gamma^{3}}{12D}\right]  \right\}  }{\sqrt{2\pi}2^{\frac{1}{3}%
}D^{\frac{2}{3}}}\frac{1}{2\pi i}\int\limits_{Br}\exp\left\{  2^{-\frac{2}{3}%
}D^{-\frac{1}{3}}\gamma\lambda\right\}  \frac{\operatorname{Ai}\left(
\lambda+2^{-\frac{1}{3}}D^{-\frac{2}{3}}\mu\right)  }{\left[
\operatorname{Ai}\left(  \lambda\right)  \right]  ^{2}}d\lambda
\]
where $Br$ is a vertical contour in the complex $\lambda$-plane on which
$\operatorname{Re}(\lambda)\geq0$ and $\operatorname{Ai}(\cdot)$ is the Airy function.

By combining the preceding results we have, on the corner scale,
\begin{align}
F(x,\eta)  &  \sim\varepsilon^{-\frac{7}{6}}\exp\left\{  \Psi_{C}(\mu
,\gamma)\right\}  L_{C}(\mu,\gamma)\equiv\widetilde{F}(\mu,\gamma)\nonumber\\
\Psi_{C}(\mu,\gamma)  &  =-\frac{\eta^{2}}{2\varepsilon}+\frac{\mu\gamma}%
{2D}-\frac{\gamma^{3}}{12D}\label{Fcorner}\\
L_{C}(\mu,\gamma)  &  =\frac{1}{\sqrt{2\pi}2^{\frac{1}{3}}D^{\frac{2}{3}}%
}\frac{1}{2\pi i}\int\limits_{Br}\exp\left\{  2^{-\frac{2}{3}}D^{-\frac{1}{3}%
}\gamma\lambda\right\}  \frac{\operatorname{Ai}\left(  \lambda+2^{-\frac{1}%
{3}}D^{-\frac{2}{3}}\mu\right)  }{\left[  \operatorname{Ai}\left(
\lambda\right)  \right]  ^{2}}d\lambda\nonumber
\end{align}
where $\mu=\varepsilon^{-\frac{2}{3}}x$ and $\gamma=\varepsilon^{-\frac{1}{3}%
}\left(  \eta-1\right)  .$

In \cite[Theorem 4]{knesslcorner} several asymptotic expansions for a function
closely related to (\ref{Fcorner}) were obtained. We use these results in the
following sections in order to match the different solutions that we have
found so far, and determine the unknown functions and constants.

\subsection{Matching the solution in Region II and the inner solution}

From (\ref{rayII1}) we get the local inversion between $(\tau,\sigma)$ and
$(x,\eta),$ for $x\rightarrow0,$%
\begin{align}
\tau &  \sim\sqrt{2}D^{\frac{1}{4}}\beta(\eta)^{-\frac{1}{4}}x^{\frac{1}{2}%
}+\frac{2}{3}\frac{\alpha(\eta)}{\beta(\eta)}x+\frac{\sqrt{2}}{36}\beta
(\eta)^{-\frac{7}{4}}D^{-\frac{1}{4}}\left[  14\beta(\eta)+11D\beta
(\eta)-20D\right]  x^{\frac{3}{2}}\label{rayIIlocal}\\
\sigma &  \sim\eta-\sqrt{2}D^{\frac{1}{4}}\beta(\eta)^{-\frac{1}{4}}%
x^{\frac{1}{2}}+\frac{1}{3}\frac{\alpha(\eta)}{\sqrt{D\beta(\eta)}%
}x\nonumber\\
&  +\frac{\sqrt{2}}{36}\beta(\eta)^{-\frac{5}{4}}D^{-\frac{3}{4}}\left[
10\beta(\eta)+D\beta(\eta)-4D\right]  x^{\frac{3}{2}}.\nonumber
\end{align}
Using (\ref{rayIIlocal}) in (\ref{L}) we obtain
\begin{equation}
L\sim L_{0}(\eta)2^{-\frac{1}{4}}\left[  \frac{D}{\beta(\eta)}\right]
^{\frac{1}{8}}x^{-\frac{1}{4}}+O(x^{\frac{1}{4}}),\quad x\rightarrow0.
\label{Llocal}%
\end{equation}

\bigskip Expanding (\ref{R0sol}) for $\mu\rightarrow\infty$ yields
\begin{align}
R_{0}  &  \sim C_{1}(\eta)\beta(\eta)^{-\frac{1}{24}}2^{-\frac{11}{12}%
}D^{\frac{5}{24}}\frac{1}{\sqrt{\pi}}\mu^{-\frac{1}{4}}\label{R0inf}\\
&  \times\exp\left[  -\frac{\sqrt{2}}{3}\beta(\eta)^{\frac{1}{4}}D^{-\frac
{5}{4}}\mu^{\frac{3}{2}}-\sqrt{2}\beta(\eta)^{\frac{1}{4}}D^{-\frac{1}{4}%
}\Gamma^{\prime}(\eta)\mu^{\frac{1}{2}}\right]  .\nonumber
\end{align}
Using (\ref{R0inf}) in (\ref{F2/3}) yields the expansion of $F$ in
(\ref{F2/3}) as $\mu\rightarrow\infty.$ By expanding $\Phi(x,\eta)$ for small
$x$ we see that the exponential parts match automatically and the matching of
the algebraic factors implies that
\[
\varepsilon^{\nu_{2}}L_{0}(\eta)2^{-\frac{1}{4}}\left[  \frac{D}{\beta(\eta
)}\right]  ^{\frac{1}{8}}x^{-\frac{1}{4}}=\varepsilon^{\nu_{3}}C_{1}%
(\eta)D^{\frac{5}{24}}\beta(\eta)^{-\frac{1}{24}}2^{-\frac{11}{12}}\frac
{1}{\sqrt{\pi}}\mu^{-\frac{1}{4}}.
\]
Hence we have
\begin{align}
\nu_{2}  &  =\nu_{3}+\frac{1}{6}\label{match1}\\
L_{0}(\eta)  &  =C_{1}(\eta)D^{\frac{1}{12}}\beta(\eta)^{\frac{1}{12}%
}2^{-\frac{2}{3}}\frac{1}{\sqrt{\pi}}.\nonumber
\end{align}

\subsection{Matching the inner and inner-inner solutions}

We take the limit $\mu\rightarrow0$ in (\ref{R0sol}) and (\ref{R1sol}) to get
\begin{equation}
R_{0}\sim C_{1}(\eta)\left\{  \operatorname{Ai}\left[  2^{\frac{2}{3}}%
D^{\frac{1}{6}}\beta\left(  \eta\right)  ^{\frac{1}{6}}\Gamma^{\prime}%
(\eta)\right]  +\frac{1}{2}2^{\frac{2}{3}}D^{-\frac{5}{6}}\beta\left(
\eta\right)  ^{\frac{1}{6}}\operatorname{Ai}^{\prime}\left[  2^{\frac{2}{3}%
}D^{\frac{1}{6}}\beta\left(  \eta\right)  ^{\frac{1}{6}}\Gamma^{\prime}%
(\eta)\right]  \mu\right\}  \label{R00}%
\end{equation}
and
\begin{align}
R_{1}  &  \sim\left[  \frac{1}{6}\sqrt{\frac{D}{\beta(\eta)}}C_{1}(\eta
)\beta^{\prime}(\eta)\left(  \Gamma^{\prime}(\eta)\right)  ^{2}+C_{2}%
(\eta)+2\sqrt{D\beta(\eta)}C_{1}(\eta)\Gamma^{\prime\prime}(\eta
)\Gamma^{\prime}(\eta)\right] \nonumber\\
&  \times\operatorname{Ai}\left[  2^{\frac{2}{3}}D^{\frac{1}{6}}\beta
(\eta)^{\frac{1}{6}}\Gamma^{\prime}(\eta)\right]  +\label{R10}\\
&  +2^{\frac{2}{3}}\left[  D^{\frac{1}{6}}\beta(\eta)^{\frac{1}{6}}%
C_{1}^{\prime}(\eta)-\frac{1}{2}D^{\frac{2}{3}}\beta(\eta)^{-\frac{1}{3}}%
C_{1}(\eta)+\frac{1}{3}D^{\frac{1}{6}}C_{1}(\eta)\alpha\beta^{-\frac{5}{6}%
}\right] \nonumber\\
&  \times\operatorname{Ai}^{\prime}\left[  2^{\frac{2}{3}}D^{\frac{1}{6}}%
\beta(\eta)^{\frac{1}{6}}\Gamma^{\prime}(\eta)\right]  .\nonumber
\end{align}

In order to complete the matching with the inner-inner solution, we must have
\[
\varepsilon^{\nu_{4}}\left.  w(v,\eta)\right|  _{v\rightarrow\infty}%
\sim\varepsilon^{\nu_{3}}\left[  R_{0}(\mu,\eta)+\varepsilon^{\frac{1}{3}%
}R_{1}(\mu,\eta)\right]  _{\mu\rightarrow0}.
\]
From (\ref{Wsol}), (\ref{R00}) and (\ref{R10}) we conclude that
\[
\nu_{3}+\frac{1}{3}=\nu_{4}%
\]%

\begin{equation}
\operatorname{Ai}\left[  2^{\frac{2}{3}}D^{\frac{1}{6}}\beta\left(
\eta\right)  ^{\frac{1}{6}}\Gamma^{\prime}(\eta)\right]  =0 \label{Airoot1}%
\end{equation}%

\begin{equation}
C_{1}(\eta)\frac{1}{2}2^{\frac{2}{3}}D^{-\frac{5}{6}}\beta\left(  \eta\right)
^{\frac{1}{6}}\operatorname{Ai}^{\prime}\left[  2^{\frac{2}{3}}D^{\frac{1}{6}%
}\beta\left(  \eta\right)  ^{\frac{1}{6}}\Gamma^{\prime}(\eta)\right]
=w(\eta)\frac{1}{2D}(\eta-1) \label{w1}%
\end{equation}%

\begin{align}
w(\eta)  &  =2^{\frac{2}{3}}\left[  D^{\frac{1}{6}}\beta\left(  \eta\right)
^{\frac{1}{6}}C_{1}^{\prime}\left(  \eta\right)  -\frac{1}{2}D^{\frac{2}{3}%
}\beta(\eta)^{-\frac{1}{3}}C_{1}(\eta)\right. \label{w2}\\
&  \left.  +\frac{1}{3}D^{\frac{1}{6}}C_{1}(\eta)\alpha\beta^{-\frac{5}{6}%
}\right]  \operatorname{Ai}^{\prime}\left[  2^{\frac{2}{3}}D^{\frac{1}{6}%
}\beta(\eta)^{\frac{1}{6}}\Gamma^{\prime}(\eta)\right]  .\nonumber
\end{align}

If we denote by $r_{0}$ the smallest (in absolute value) of the roots of
$\operatorname{Ai},$ i.e.,
\[
r_{0}=\max\left\{  z:\operatorname{Ai}(z)=0\right\}  \simeq-2.33810741
\]
then we have from (\ref{Airoot1})
\begin{equation}
\Gamma^{\prime}(\eta)=2^{-\frac{2}{3}}D^{-\frac{1}{6}}\beta\left(
\eta\right)  ^{-\frac{1}{6}}r_{0}. \label{Gammaderiv}%
\end{equation}
From (\ref{w1}) and (\ref{w2}) we obtain an ODE for $C_{1}(\eta)$%
\begin{equation}
C_{1}^{\prime}(\eta)+\left[  \frac{1}{4}\frac{\alpha\left(  \eta\right)
}{D\beta\left(  \eta\right)  }-\frac{1}{4\sqrt{D\beta\left(  \eta\right)  }%
}-\frac{1}{6}\frac{\alpha\left(  \eta\right)  }{\beta(\eta)}-\frac{1}{\eta
-1}\right]  C_{1}(\eta)=0, \label{C1ode}%
\end{equation}
and a relation between $C_{1}(\eta)$ and $w(\eta)$%
\begin{equation}
w(\eta)=\frac{1}{\eta-1}C_{1}(\eta)2^{\frac{2}{3}}D^{\frac{1}{6}}\beta\left(
\eta\right)  ^{\frac{1}{6}}\operatorname{Ai}^{\prime}(r_{0}),\quad\eta>1.
\label{C1w}%
\end{equation}
The solution of (\ref{C1ode}) is
\begin{equation}
C_{1}(\eta)=k_{0}(\eta-1)\beta\left(  \eta\right)  ^{-\frac{1}{6}}\left[
\frac{\alpha\left(  \eta\right)  }{\sqrt{D+1}}+\sqrt{\beta\left(  \eta\right)
}\right]  ^{\frac{\sqrt{D}}{2\sqrt{D+1}}} \label{C1sol}%
\end{equation}
with $k_{0}$ a constant to be determined.

\subsection{Matching the corner and Region I solutions}

From \cite[Theorem 4 (i)]{knesslcorner} we have the following result valid
when $\mu$ and/or $\left|  \gamma\right|  \rightarrow\infty$ with
$\gamma-\sqrt{\mu}\rightarrow-\infty,$
\begin{align}
\widetilde{F}(\mu,\gamma)  &  \sim\varepsilon^{-\frac{7}{6}}L_{I}(\mu
,\gamma)\exp\left\{  \Psi_{I}(\mu,\gamma)\right\} \label{cornerandI}\\
\Psi_{I}(\mu,\gamma)  &  =-\frac{1}{27D}\left\{  \gamma^{3}-18\mu
\gamma+\left[  \gamma^{2}+6\mu\right]  ^{\frac{3}{2}}\right\}  -\frac
{1}{2\varepsilon}\eta^{2}\\
L_{I}(\mu,\gamma)  &  =\frac{1}{D}\frac{1}{\sqrt{\pi}}\frac{\sqrt{6}}%
{18}\left[  \gamma^{2}+6\mu\right]  ^{-\frac{1}{4}}\left[  \sqrt{\gamma
^{2}+6\mu}-2\mu\right]  ^{\frac{3}{2}}.\nonumber
\end{align}

We can invert the ray transformation (\ref{ray3}) locally when $x\rightarrow
0,\ \eta\rightarrow1$ to get
\begin{equation}
t\sim\frac{1}{3}\left[  z+(\eta-1)\right]  ,\quad s\sim1-\frac{1}{3}z+\frac
{2}{3}(\eta-1) \label{ts}%
\end{equation}
where
\[
z=\sqrt{(\eta-1)^{2}+6x}.
\]
Using (\ref{ts}) in the ray expansion (\ref{prop2}) yields, as $(x,\eta
)\rightarrow(0,1),$
\begin{align*}
F  &  \sim\varepsilon^{-\frac{3}{2}}K(x,\eta)\exp\left[  \frac{1}{\varepsilon
}\Psi(x,\eta)\right] \\
K  &  \sim\frac{1}{D}\frac{1}{\sqrt{2\pi}}\frac{1}{\sqrt{z}}\left\{  \frac
{1}{3}\left[  z-2(\eta-1)\right]  \right\}  ^{\frac{3}{2}}\\
\Psi+\frac{1}{2}\eta^{2}  &  \sim\frac{1}{D}\left\{  -\frac{1}{27}z^{3}%
+\frac{1}{9}(\eta-1)z^{2}-\frac{4}{27}(\eta-1)^{3}\right\}  ,
\end{align*}
which agrees with (\ref{cornerandI}).

\subsection{Matching the corner and Region II solutions}

From \cite[Theorem 4 (iv) ]{knesslcorner} we have
\begin{align}
\widetilde{F}(\mu,\gamma)  &  \sim\varepsilon^{-\frac{7}{6}}L_{II}(\mu
,\gamma)\exp\left\{  \Psi_{II}(\mu,\gamma)\right\} \label{LcII}\\
L_{II}(\mu,\gamma)  &  =D^{-\frac{5}{6}}\frac{1}{\pi}\frac{1}{\left[
\operatorname{Ai}^{\prime}(r_{0})\right]  ^{2}}2^{-\frac{29}{12}}\gamma
\mu^{-\frac{1}{4}}%
\end{align}%

\begin{align}
\Psi_{II}(\mu,\gamma)  &  =-\frac{1}{2\varepsilon}\eta^{2}-\frac{1}{12D}%
\gamma^{3}+\frac{1}{2D}\mu\gamma-\frac{1}{3D}\sqrt{2}\mu^{\frac{3}{2}%
}\nonumber\\
&  +\frac{1}{2}2^{\frac{1}{3}}D^{-\frac{1}{3}}r_{0}\gamma-2^{-\frac{1}{6}%
}D^{-\frac{1}{3}}r_{0}\sqrt{\mu} \label{psiCII}%
\end{align}
which is valid when $\mu$ and $\gamma\rightarrow\infty,$ with $\gamma
-\sqrt{\mu}\rightarrow\infty.$

Combining (\ref{Llocal}), (\ref{match1}) and (\ref{C1sol}) we have
\[
L\sim k_{0}2^{-\frac{11}{12}}\frac{1}{\sqrt{\pi}}\left[  \frac{D}{\sqrt{D+1}%
}+\sqrt{D}\right]  ^{\frac{\sqrt{D}}{2\sqrt{D+1}}}(\eta-1)x^{-\frac{1}{4}%
},\quad x\rightarrow0,
\]
which agrees with (\ref{LcII}) if
\begin{equation}
k_{0}=D^{-\frac{5}{6}}\frac{1}{\sqrt{\pi}}\frac{1}{\left[  \operatorname{Ai}%
^{\prime}(r_{0})\right]  ^{2}}2^{-\frac{3}{2}}\left[  \frac{D}{\sqrt{D+1}%
}+\sqrt{D}\right]  ^{-\frac{\sqrt{D}}{2\sqrt{D+1}}}. \label{k0}%
\end{equation}

Since in Region II $F(x,\eta)\sim\varepsilon^{\nu_{2}}\exp\left[
\varepsilon^{-1}\Phi(x,\eta)+\varepsilon^{-\frac{1}{3}}\Gamma(x,\eta)\right]
L(x,\eta),$ we must have
\begin{equation}
\nu_{2}=-\frac{4}{3}. \label{nu2}%
\end{equation}
We use (\ref{rayIIlocal}) in (\ref{rayII1}), (\ref{phi00}) and
(\ref{Gammaderiv}) and find that, as $(x,\eta)\rightarrow(0,1),$
\begin{align*}
\Phi(x,\eta)  &  \sim-\frac{1}{2}-(\eta-1)-\frac{1}{2}(\eta-1)^{2}-\frac
{1}{12D}(\eta-1)^{3}+\frac{1}{2D}x(\eta-1)-\frac{1}{3D}\sqrt{2}x^{\frac{3}{2}%
}\\
\Gamma(\sigma)  &  \sim\Gamma(1)+\frac{1}{2}2^{\frac{1}{3}}D^{-\frac{1}{3}%
}r_{0}(\eta-1)-2^{-\frac{1}{6}}D^{-\frac{1}{3}}r_{0}\sqrt{x}%
\end{align*}
and from (\ref{psiCII}) we conclude that
\begin{equation}
\Gamma(1)=0. \label{gamma1}%
\end{equation}

We have now determined all the unknown functions from the previous sections
and these we summarize below%

\begin{equation}
L(x,\eta)=D^{-\frac{3}{4}}(\sigma-1)\frac{1}{\pi}2^{-\frac{5}{2}}\beta\left(
\sigma\right)  ^{-\frac{1}{12}}\left[  \frac{\alpha\left(  \sigma\right)
+\sqrt{\beta\left(  \sigma\right)  (D+1)}}{D+\sqrt{D(D+1)}}\right]
^{\frac{\sqrt{D}}{2\sqrt{D+1}}}\frac{1}{\left[  \operatorname{Ai}^{\prime
}(r_{0})\right]  ^{2}}\frac{e^{\frac{\tau}{2}}}{\sqrt{\widetilde{J}}}
\label{Ldef}%
\end{equation}%

\begin{equation}
\Gamma(\sigma)=2^{-\frac{2}{3}}D^{-\frac{1}{6}}r_{0}\int\limits_{1}^{\sigma
}\beta\left(  u\right)  ^{-\frac{1}{6}}du \label{Gammadef}%
\end{equation}%

\begin{align}
R_{0}(\mu,\eta)  &  =(\eta-1)D^{-\frac{5}{6}}\frac{1}{\sqrt{\pi}}2^{-\frac
{3}{2}}\beta(\eta)^{-\frac{1}{6}}\left[  \frac{\alpha(\eta)+\sqrt{\beta
(\eta)(D+1)}}{D+\sqrt{D(D+1)}}\right]  ^{\frac{\sqrt{D}}{2\sqrt{D+1}}%
}\label{R0def}\\
\times &  \frac{\operatorname{Ai}\left[  2^{-\frac{1}{3}}D^{-\frac{5}{6}}%
\beta(\eta)^{\frac{1}{6}}\mu+r_{0}\right]  }{\left[  \operatorname{Ai}%
^{\prime}(r_{0})\right]  ^{2}}\nonumber
\end{align}%

\begin{equation}
W(v,\eta)=2^{-\frac{5}{6}}\frac{1}{\sqrt{\pi}}D^{-\frac{2}{3}}\left[
\frac{\alpha(\eta)+\sqrt{\beta(\eta)(D+1)}}{D+\sqrt{D(D+1)}}\right]
^{\frac{\sqrt{D}}{2\sqrt{D+1}}}\frac{1}{\operatorname{Ai}^{\prime}(r_{0}%
)}\left[  \frac{1}{2D}(\eta-1)v+1\right]  . \label{Wdef}%
\end{equation}

With (\ref{Ldef}) and (\ref{Gammadef}) we have completely determined the ray
expansion in Region II, with (\ref{R0def}) we have the inner solution (for
$x=O(\varepsilon^{\frac{2}{3}})$ and $\eta>1$) and with (\ref{Wdef}) we have
the inner-inner solution (for $x=O(\varepsilon)$ and $\eta>1$). We have also
show that
\[
\nu_{2}=-\frac{4}{3},\ \nu_{3}=-\frac{3}{2}\text{ and }\nu_{4}=-\frac{7}{6}.
\]

\section{Transition layer}

Finally we shall find the boundary layer solution near the curve $x=X_{0}%
(\eta)$ defined by (\ref{specialrays}), which separates Regions I and II. We
introduce the stretched variable $\omega=(x-X_{0})\varepsilon^{-\frac{1}{3}}$
and (\ref{scaled}) becomes
\begin{equation}
-2\eta^{2}(\eta-1)F_{\omega}+\eta^{2}(\eta F_{\eta}+F)\varepsilon^{\frac{1}%
{3}}+\beta F_{\omega\omega}\varepsilon^{\frac{2}{3}}-\left[  2\eta
(\eta-1)F_{\omega\eta}+F_{\omega}\right]  \varepsilon+\eta^{2}F_{\eta\eta
}\varepsilon^{\frac{4}{3}}=0. \label{corner1}%
\end{equation}

When $s=1\ (\sigma=1),$ $t=\ln(\eta)\ \left(  \tau=\ln(\eta)\right)  $ and we
have
\begin{equation}
j=J\left[  \ln(\eta),1\right]  =2\left(  1+\frac{1}{D}\right)  \ln\left(
\eta\right)  \eta+\frac{1}{D}\left(  4-3\eta-\frac{1}{\eta}\right)
=2\widetilde{J}\left[  \ln(\eta),1\right]  =2j_{1}. \label{j}%
\end{equation}
Since
\[
\Psi\sim-\frac{1}{2}\eta^{2}-\frac{\eta}{2Dj}(x-X_{0})^{2},\quad x\rightarrow
X_{0}%
\]
we should look for solutions of the form
\begin{equation}
F\sim\varepsilon^{\nu_{6}}\exp\left\{  -\frac{1}{2\varepsilon}\eta^{2}%
-\frac{\eta}{2Dj}\omega^{2}\varepsilon^{-\frac{1}{3}}\right\}  \Upsilon
(\omega,\eta). \label{corner2}%
\end{equation}
Using (\ref{corner2}) in (\ref{corner1}) yields for $\Upsilon$ the equation
\[
2D^{2}j^{2}\omega\beta(\eta)\Upsilon_{w}+\eta^{2}D^{3}j^{3}\Upsilon_{\eta
}+\beta(\eta)\left[  D^{2}j^{2}-2\omega^{3}(\eta-1)\right]  \Upsilon=0
\]
whose general solution is
\begin{equation}
\Upsilon(\omega,\eta)=g\left(  \frac{\eta\omega}{Dj}\right)  \sqrt{\frac{\eta
}{Dj}}\exp\left\{  -\frac{\omega^{3}}{2\eta D^{3}j^{3}}\left[  \left(
2\eta-1\right)  \left(  2D\eta^{2}+2\eta^{2}-2\eta+1\right)  \right]
\right\}  , \label{Y}%
\end{equation}
where $g$ is a function still unknown. It will be determined in the next
section by matching with the corner solution.

\subsection{Matching the corner and transition layer solutions}

Let us first introduce the new variable $\Omega$ defined by
\begin{equation}
\Omega=\frac{1}{(2D)^{\frac{1}{3}}}\left(  \mu-\frac{1}{2}\gamma^{2}\right)
\frac{1}{\gamma}. \label{bigomega}%
\end{equation}

From \cite[Theorem 4 (ii) ]{knesslcorner} we have the following result, for
$\mu,\gamma\rightarrow\infty,$ $\Omega\ $fixed
\begin{equation}
\varepsilon^{-\frac{7}{6}}e^{\frac{\eta^{2}}{2\varepsilon}}F(\mu,\gamma
)\sim\varepsilon^{-\frac{7}{6}}\frac{2^{\frac{5}{6}}}{4\pi\sqrt{D\gamma}}%
\wp(\Omega)\exp\left\{  \frac{\Omega^{3}}{6}-\frac{1}{4}\gamma\Omega
^{2}2^{\frac{2}{3}}D^{-\frac{1}{3}}\right\}  \label{cornertran}%
\end{equation}
where
\[
\wp(\Omega)=\frac{1}{2\pi i}\int\limits_{Br}\frac{e^{-\lambda\Omega}}{\left[
\operatorname{Ai}\left(  2^{\frac{1}{3}}\lambda\right)  \right]  ^{2}}%
d\lambda.
\]
The following properties of $\wp(\Omega)$ are established in
\cite{knesslcorner} \
\begin{align}
\wp(0)  &  =2^{-\frac{1}{3}}\label{P}\\
\wp(\Omega)  &  \sim\Omega^{\frac{3}{2}}\sqrt{\pi}2^{-\frac{5}{6}}\exp\left\{
-\frac{\Omega^{3}}{24}\right\}  ,\quad\Omega\rightarrow\infty\nonumber\\
\wp(\Omega)  &  \sim-\frac{\Omega2^{-\frac{2}{3}}}{\left[  \operatorname{Ai}%
^{\prime}(r_{0})\right]  ^{2}}\exp\left\{  -2^{-\frac{1}{3}}r_{0}%
\Omega\right\}  ,\quad\Omega\rightarrow-\infty.\nonumber
\end{align}

In order to match with (\ref{Y}), we first note that
\begin{equation}
\omega\sim(2D)^{\frac{1}{3}}(\eta-1)\Omega,\quad\eta\rightarrow1 \label{omega}%
\end{equation}
thus the right side of (\ref{corner2}) behaves as
\begin{align*}
&  \varepsilon^{\nu_{6}}g\left[  \left(  2D\right)  ^{-\frac{2}{3}}%
\Omega\right]  \frac{1}{\sqrt{2D(\eta-1)}}\exp\left\{  -\frac{1}{8}\frac
{2D+1}{D^{2}}\Omega^{3}-\frac{1}{4}\Omega^{2}2^{\frac{2}{3}}D^{-\frac{1}{3}%
}(\eta-1)\varepsilon^{-\frac{1}{3}}\right\} \\
&  =\varepsilon^{\nu_{6}}g\left[  \left(  2D\right)  ^{-\frac{2}{3}}%
\Omega\right]  \frac{1}{\varepsilon^{\frac{1}{6}}\sqrt{2D\gamma}}\exp\left\{
-\frac{1}{8}\frac{2D+1}{D^{2}}\Omega^{3}-\frac{1}{4}\Omega^{2}2^{\frac{2}{3}%
}D^{-\frac{1}{3}}\gamma\right\}  .
\end{align*}

Comparing the above with (\ref{cornertran}) we must have
\[
\nu_{6}=-1
\]
and
\[
g\left[  \left(  2D\right)  ^{-\frac{2}{3}}\Omega\right]  =\exp\left\{
\frac{\Omega^{3}}{6}+\frac{1}{8}\frac{2D+1}{D^{2}}\Omega^{3}\right\}  \frac
{1}{\pi}2^{-\frac{2}{3}}\wp(\Omega)
\]
which implies that
\[
g(Z)=\exp\left\{  \frac{Z^{3}}{6}\left(  4D^{2}+6D+3\right)  \right\}
\frac{1}{\pi}2^{-\frac{2}{3}}\wp\left[  \left(  2D\right)  ^{\frac{2}{3}%
}Z\right]  .
\]

We conclude by writing the complete transition layer solution in
(\ref{corner2})%

\begin{align}
F  &  \sim\frac{1}{\varepsilon\pi}2^{-\frac{2}{3}}\sqrt{\frac{\eta}{Dj}}%
\wp\left[  \frac{2^{\frac{2}{3}}}{D^{\frac{1}{3}}}\frac{\eta\omega}{j}\right]
\exp\left\{  -\frac{\eta^{2}}{2\varepsilon}-\frac{\eta}{2Dj\varepsilon
^{\frac{1}{3}}}\omega^{2}+\frac{1}{6}\left(  4D^{2}+6D+3\right)  \left(
\frac{\eta\omega}{Dj}\right)  ^{3}\right. \nonumber\\
&  \left.  -\frac{\omega^{3}}{2\eta D^{3}j^{3}}\left[  \left(  2\eta-1\right)
\left(  2D\eta^{2}+2\eta^{2}-2\eta+1\right)  \right]  \right\}
\label{cornerlayersol}\\
&  \equiv\varepsilon^{-1}L_{X_{0}}(\omega,\eta)\exp\left\{  \Psi_{X_{0}%
}(\omega,\eta;\varepsilon)\right\}  .\nonumber
\end{align}
In the next two subsections we will show that (\ref{cornerlayersol}) matches
to both of the solutions in Regions I and II.

\subsection{Matching the solution in Region I and the transition layer solution}

When $x$ is close to $X_{0}$, we can invert the equations (\ref{ray3}) to get
\begin{equation}
t\sim\ln(\eta),\quad s\sim1-\frac{\eta}{j}(x-X_{0}). \label{tranIlocal}%
\end{equation}
Using the above in (\ref{prop2}) we obtain
\begin{align}
\frac{1}{\varepsilon}\Psi &  \sim\frac{1}{\varepsilon}\left[  -\frac{1}{2}%
\eta^{2}-\frac{\eta}{2Dj}(x-X_{0})^{2}+\frac{1}{2\eta D^{3}j^{3}}\beta
(\eta)^{2}(x-X_{0})^{3}\right]  ,\label{transitionandI}\\
\varepsilon^{-\frac{3}{2}}K  &  \sim\varepsilon^{-\frac{3}{2}}\frac{\eta^{2}%
}{D\sqrt{2\pi}j^{2}}\left(  x-X_{0}\right)  ^{\frac{3}{2}}.\nonumber
\end{align}

From the definition of $\Omega$ in (\ref{bigomega}) we see that
\[
\Omega=(2D)^{-\frac{1}{3}}\left[  \frac{1}{\eta-1}\left(  \omega
+\varepsilon^{-\frac{1}{3}}X_{0}\right)  -\frac{1}{2}\varepsilon^{-\frac{1}%
{3}}\left(  \eta-1\right)  \right]
\]
and thus $\Omega\rightarrow\pm\infty$ when $\omega\rightarrow\pm\infty.$
Expanding (\ref{cornerlayersol}) for $\omega\rightarrow\infty$ and taking into
account (\ref{P}) we have
\[
L_{_{X_{0}}}\sim\frac{1}{\sqrt{2\pi}D}\left(  \frac{\eta}{j}\right)
^{2}\omega^{\frac{3}{2}}.
\]
This agrees with (\ref{transitionandI}), since $x-X_{0}=\varepsilon^{\frac
{1}{3}}\omega$.

\subsection{Matching the solution in Region II and the transition layer solution}

For $x\rightarrow X_{0}$ we get from (\ref{rayII1})
\[
\tau\sim\ln(\eta),\quad\sigma\sim1-\frac{\eta}{j_{1}}(x-X_{0})
\]
which when used in (\ref{rayII1}), (\ref{Gammadef}) and (\ref{Ldef}) yields
\begin{align}
\frac{1}{\varepsilon}\Phi &  \sim-\frac{1}{2\varepsilon}\eta^{2}-\frac{\eta
}{4D\varepsilon j_{1}}(x-X_{0})^{2}+\frac{1}{8\eta D^{3}j_{1}^{3}}\left[
\left(  4D^{2}+6D+3\right)  \eta^{4}\right. \nonumber\\
&  \left.  -3\left(  2\eta-1\right)  \left(  2D\eta^{2}+2\eta^{2}%
-2\eta+1\right)  \right]  \frac{1}{6\varepsilon}(x-X_{0})^{3},\nonumber\\
\varepsilon^{-\frac{1}{3}}\Gamma &  \sim-\frac{1}{2}\varepsilon^{-\frac{1}{3}%
}2^{\frac{1}{3}}D^{-\frac{1}{3}}r_{0}\frac{\eta}{j_{1}}\left(  x-X_{0}\right)
,\label{tranandII}\\
\varepsilon^{-\frac{4}{3}}L  &  \sim-\varepsilon^{-\frac{4}{3}}D^{-\frac{5}%
{6}}\frac{1}{\pi}\frac{1}{\left[  \operatorname{Ai}^{\prime}(r_{0})\right]
^{2}}2^{-\frac{13}{6}}\left(  \frac{\eta}{j_{1}}\right)  ^{\frac{3}{2}}\left(
x-X_{0}\right)  .\nonumber
\end{align}
From (\ref{cornerlayersol}) and (\ref{P}) we find that when $\omega
\rightarrow-\infty$%
\[
L_{_{X_{0}}}\sim-\frac{1}{\pi}2^{-\frac{2}{3}}D^{-\frac{5}{6}}\left(
\frac{\eta}{j}\right)  ^{\frac{3}{2}}\omega\left[  \operatorname{Ai}^{\prime
}(r_{0})\right]  ^{-2}\exp\left\{  -2^{\frac{1}{3}}D^{-\frac{1}{3}}r_{0}%
\frac{\eta}{j}\omega\right\}  .
\]
This matches with (\ref{tranandII}) if we take into account (\ref{j}).

\section{Summary of results and discussion}

Below we summarize the main results of this section, which consist of the
asymptotic expansions of $F(x,\eta)$ in (\ref{scaled}) in the various parts of
the $(x,\eta)$ plane.

(A) Region I $\left\{  x>X_{0}=\eta-\ln(\eta)-1,\quad\eta>1\right\}
\cup\left\{  x>0,\quad\eta\leq1\right\}  $%

\[
F(x,\eta)\sim\varepsilon^{-\frac{3}{2}}K(x,\eta)\exp\left[  \frac
{1}{\varepsilon}\Psi(x,\eta)\right]
\]%

\[
x=e^{t}-1-t-\frac{(D+1)(2t-e^{t})+D+e^{-t}}{D}(s-1),
\]%

\[
\eta=e^{t}+\frac{e^{-t}+(D+1)e^{t}-2}{D}(s-1),
\]%

\[
K(x,\eta)=\frac{1}{\sqrt{2\pi}}(1-s)^{\frac{3}{2}}\frac{e^{\frac{t}{2}}}%
{\sqrt{J}},
\]%

\begin{align*}
J  &  =\left[  2(t-2)(s-1)D^{-2}+(-2t-5s+4ts+2)D^{-1}-s+2ts+1\right]  e^{t}\\
&  +\left[  -2(t+2)(s-1)D^{-2}+(2t-2ts+2-3s)D^{-1}\right]  e^{-t}\\
&  +8(s-1)D^{-2}+4(2s-1)D^{-1},
\end{align*}%

\begin{align*}
\Psi(x,\eta)  &  =-\frac{1}{2}e^{2t}+\frac{2e^{t}-(D+1)e^{2t}-1}{D}(s-1)\\
&  +\frac{-1+[4e^{t}-2(t+1)](D+1)-e^{2t}(D+1)^{2}}{2D^{2}}(s-1)^{2}.
\end{align*}

(B) Corner layer $x=\mu\varepsilon^{\frac{2}{3}},\quad\eta-1=\gamma
\varepsilon^{\frac{1}{3}}$%

\[
F(x,\eta)\sim\varepsilon^{-\frac{7}{6}}\exp\left\{  \Psi_{C}(\mu
,\gamma)\right\}  L_{C}(\mu,\gamma)
\]%

\[
L_{C}(\mu,\gamma)=\frac{1}{\sqrt{2\pi}\left(  2D^{2}\right)  ^{\frac{1}{3}}%
}\frac{1}{2\pi i}\int\limits_{Br}\exp\left\{  (4D)^{-\frac{1}{3}}\gamma
\lambda\right\}  \frac{\operatorname{Ai}\left[  \lambda+\left(  2D^{2}\right)
^{-\frac{1}{3}}\mu\right]  }{\left[  \operatorname{Ai}\left(  \lambda\right)
\right]  ^{2}}d\lambda,
\]%

\[
\Psi_{C}(\mu,\gamma)=-\frac{\eta^{2}}{2\varepsilon}+\frac{\mu\gamma}{2D}%
-\frac{\gamma^{3}}{12D}.
\]

(C) Transition layer $x-X_{0}=\omega\varepsilon^{\frac{1}{3}}$%

\[
F(x,\eta)\sim\varepsilon^{-1}L_{X_{0}}(\omega,\eta)\exp\left\{  \Psi_{X_{0}%
}(\omega,\eta;\varepsilon)\right\}  ,
\]%

\[
L_{X_{0}}(x,\eta)=\frac{1}{\pi}2^{-\frac{2}{3}}\sqrt{\frac{\eta}{Dj}}%
\wp\left[  \frac{2^{\frac{2}{3}}}{D^{\frac{1}{3}}}\frac{\eta\omega}{j}\right]
,
\]%

\[
j=2\left(  1+\frac{1}{D}\right)  \ln\left(  \eta\right)  \eta+\frac{1}%
{D}\left(  4-3\eta-\frac{1}{\eta}\right)  ,
\]%

\[
\wp(\Omega)=\frac{1}{2\pi i}\int\limits_{Br}\frac{e^{-\lambda\Omega}}{\left[
\operatorname{Ai}\left(  2^{\frac{1}{3}}\lambda\right)  \right]  ^{2}}%
d\lambda,
\]%

\begin{align*}
\Psi_{X_{0}}(\omega,\eta;\varepsilon)  &  =-\frac{\eta^{2}}{2\varepsilon
}-\frac{\eta}{2Dj\varepsilon^{\frac{1}{3}}}\omega^{2}+\frac{1}{6}\left(
4D^{2}+6D+3\right)  \left(  \frac{\eta\omega}{Dj}\right)  ^{3}\\
&  -\frac{\omega^{3}}{2\eta D^{3}j^{3}}\left[  \left(  2\eta-1\right)  \left(
2D\eta^{2}+2\eta^{2}-2\eta+1\right)  \right]  .
\end{align*}

(D) Region II $\left\{  0<x<X_{0}=\eta-\ln(\eta)-1,\quad\eta>1\right\}  $%

\[
F(x,\eta)\sim\varepsilon^{-\frac{4}{3}}\exp\left[  \varepsilon^{-1}\Phi
(x,\eta)+\varepsilon^{-\frac{1}{3}}\Gamma(x,\eta)\right]  L(x,\eta),
\]%

\[
a(\sigma)=\frac{1-\sigma}{2D},\quad b(\sigma)=\frac{\sigma}{2}+\frac
{\sqrt{D\sigma^{2}+\left(  \sigma-1\right)  ^{2}}}{2\sqrt{D}},
\]%

\[
x=(b-a)e^{\tau}+(a+b-\sigma)e^{-\tau}+[2a(D+1)-1]\tau-2b+\sigma,
\]%

\[
\eta=(b-a)e^{\tau}-(a+b-\sigma)e^{-\tau}+2a,
\]%

\[
\Phi=-a^{2}(D+1)\tau+2a\left(  a-b\right)  \left(  e^{\tau}-1\right)
-\frac{1}{2}\left(  a-b\right)  ^{2}\left(  e^{2\tau}-1\right)  +\Phi
_{0}(\sigma),
\]%

\[
\Phi_{0}(\sigma)=-\frac{1}{2}-\int\limits_{1}^{\sigma}b\left(  u\right)  du,
\]%

\[
\alpha(\sigma)=\left(  D+1\right)  \sigma-1,\quad\beta(\sigma)=D\sigma
^{2}+\left(  \sigma-1\right)  ^{2},
\]%

\[
r_{0}=\max\left\{  z:\operatorname{Ai}(z)=0\right\}  \simeq-2.33810741,
\]%
\[
\Gamma(\sigma)=2^{-\frac{2}{3}}D^{-\frac{1}{6}}r_{0}\int\limits_{1}^{\sigma
}\beta\left(  u\right)  ^{-\frac{1}{6}}du,
\]%

\begin{align*}
L(x,\eta)  &  =D^{-\frac{3}{4}}(\sigma-1)\frac{1}{\pi}2^{-\frac{5}{2}}%
\beta\left(  \sigma\right)  ^{-\frac{1}{12}}\left[  \frac{\alpha\left(
\sigma\right)  +\sqrt{\beta\left(  \sigma\right)  (D+1)}}{D+\sqrt{D(D+1)}%
}\right]  ^{\frac{\sqrt{D}}{2\sqrt{D+1}}}\\
&  \times\frac{1}{\left[  \operatorname{Ai}^{\prime}(r_{0})\right]  ^{2}}%
\frac{e^{\frac{\tau}{2}}}{\sqrt{\widetilde{J}}},
\end{align*}%

\begin{align*}
\widetilde{J}  &  =\left\{  \left[  -\sigma+1+\frac{1}{2}\tau(\sigma
-1)\right]  D^{-2}+\left[  \frac{1}{2}\sqrt{\beta(\sigma)}(\tau-1)\right]
D^{-\frac{3}{2}}\right. \\
&  \left.  +\left(  -\sigma-\frac{1}{2}\tau+\tau\sigma\right)  D^{-1}+\frac
{1}{2}\tau\sqrt{\beta(\sigma)}D^{-\frac{1}{2}}+\frac{1}{2}\tau\sigma\right\}
e^{\tau}\\
&  +\left\{  \left(  \frac{1}{2}\tau+1\right)  (1-\sigma)D^{-2}+\left[
\frac{1}{2}\sqrt{\beta(\sigma)}(\tau+1)\right]  D^{-\frac{3}{2}}\right. \\
&  \left.  +\left(  -\sigma+\frac{1}{2}\tau-\tau\sigma\right)  D^{-1}+\frac
{1}{2}\tau\sqrt{\beta(\sigma)}D^{-\frac{1}{2}}-\frac{1}{2}\tau\sigma\right\}
e^{-\tau}\\
&  +2(\sigma-1)D^{-2}+2\sigma D^{-1}.
\end{align*}

(E) Inner layer $x=\mu\varepsilon^{\frac{2}{3}},\quad\eta>1$%

\[
F(x,\eta)\sim\varepsilon^{-\frac{3}{2}}\exp\left\{  \varepsilon^{-1}\left[
\Phi_{0}(\eta)+\frac{\left(  \eta-1\right)  }{2D}x\right]  +\varepsilon
^{-\frac{1}{3}}\Gamma\left(  \eta\right)  \right\}  R_{0}(\mu,\eta),
\]%

\begin{align*}
R_{0}(\mu,\eta)  &  =(\eta-1)D^{-\frac{5}{6}}\frac{1}{\sqrt{\pi}}2^{-\frac
{3}{2}}\beta(\eta)^{-\frac{1}{6}}\left[  \frac{\alpha(\eta)+\sqrt{\beta
(\eta)(D+1)}}{D+\sqrt{D(D+1)}}\right]  ^{\frac{\sqrt{D}}{2\sqrt{D+1}}}\\
&  \times\frac{\operatorname{Ai}\left[  2^{-\frac{1}{3}}D^{-\frac{5}{6}}%
\beta(\eta)^{\frac{1}{6}}\mu+r_{0}\right]  }{\left[  \operatorname{Ai}%
^{\prime}(r_{0})\right]  ^{2}}.
\end{align*}

(F) Inner-inner layer $x=v\varepsilon,\quad\eta>1$%
\[
F(x,\eta)\sim\varepsilon^{-\frac{7}{6}}\exp\left\{  \varepsilon^{-1}\left[
\Phi_{0}(\eta)+\frac{1}{2}\frac{\eta-1}{D}x\right]  +\varepsilon^{-\frac{1}%
{3}}\Gamma(\eta)\right\}  W(v,\eta),
\]%
\begin{align*}
W(v,\eta)  &  =2^{-\frac{5}{6}}\frac{1}{\sqrt{\pi}}D^{-\frac{2}{3}}\left[
\frac{\alpha(\eta)+\sqrt{\beta(\eta)(D+1)}}{D+\sqrt{D(D+1)}}\right]
^{\frac{\sqrt{D}}{2\sqrt{D+1}}}\frac{1}{\operatorname{Ai}^{\prime}(r_{0})}\\
&  \times\left[  \frac{1}{2D}(\eta-1)v+1\right]  .
\end{align*}

In that part of Region I outside the caustic region (cf. Figure 3.2) the
mapping between $(t,s)$ and $(x,\eta)$ is one-to-one, and $K$ and $\Psi$ are
unambiguously determined by the formulas in (A).

Inside the caustic region the mapping is three-to-one and we should re-write
(\ref{Ray1}) as
\[
\varepsilon^{-\frac{3}{2}}\left[  K_{1}\exp\left(  \frac{1}{\varepsilon}%
\Psi_{1}\right)  +K_{2}\exp\left(  \frac{1}{\varepsilon}\Psi_{2}\right)
+K_{3}\exp\left(  \frac{1}{\varepsilon}\Psi_{3}\right)  \right]
\]
where $\Psi_{j}$ and $K_{j\text{ }}$correspond to the three different values
of $(t,s)$ leading to the same $(x,\eta).$ When $t=0$ let us define the
starting points on the $\eta$-axis of these three rays by the ordering
$s_{1}<s_{2}<s_{3},$ where $s_{j}$ corresponds to $\Psi_{j}$ and $K_{j}.$ We
denote the two caustics by $C_{+}$ and $C_{-}$ and the cusp where they meet as
$(x_{c},\eta_{c})$. Note that the cusp location depends only on $D$.

The curve $C_{+}$ has $\eta\rightarrow-\infty$ as $x\rightarrow\infty$, while
$C_{-}$ reaches the $\eta$-axis at some critical point $(0,\eta_{\ast})$ where
again $\eta_{\ast}=\eta_{\ast}(D)$. We have verified numerically that along
$C_{+}$ we have $s_{1}=s_{2},$ $\Psi_{1}=\Psi_{2}$ and $K_{1},\ K_{2}$ develop
singularities. However, here $\Psi_{3}>$ $\Psi_{1}=\Psi_{2}$ and $K_{3}$
remains finite. Thus we have $F\sim\varepsilon^{-\frac{3}{2}}K_{1}\exp\left(
\frac{1}{\varepsilon}\Psi_{1}\right)  $ on and near $C_{+}.$ Similarly, along
$C_{-}$ we have $s_{2}=s_{3},$ $\Psi_{2}=\Psi_{3}$ and $K_{2},\ K_{3}$ develop
singularities. But $\Psi_{1}>\Psi_{2}=\Psi_{3}$ and $K_{1}$ remains finite.
Thus the result in (A) remains valid near the caustics, except near the cusp
point where all three $\Psi_{j}$ are approximately equal. Here the expansion
in (A) breaks down.

Our preliminary results suggest that a new expansion must be constructed near
the cusp with the scaling
\[
x-x_{c}=O(\sqrt{\varepsilon}),\quad\eta-\eta_{c}-A_{c}(x-x_{c})=O(\varepsilon
^{\frac{3}{4}}).
\]
Here $A_{c}$ is the slope at which both $C_{+}$ and $C_{-}$ hit the cusp. We
have thus far not been able to complete this analysis. We also note that while
the expansion near the cusp presents an interesting problem in asymptotics, it
is not needed for computing the marginal distribution $M(x)$ (\ref{mm}), which
is the most important quantity from the point of view of applications, and
which we calculate in the next section.

\section{Marginal distributions}

The last ``piece of the puzzle'', is to verify that (\ref{marginal}) is
satisfied, and also to compute the marginal distribution $M(x)$ in (\ref{mm}).

We evaluate the integral in (\ref{marginal}) for $\varepsilon\rightarrow0$.
For $\eta<1$, $F(x,\eta)$ is concentrated near $x=0$, and the result follows
from the approximation (\ref{prop1}). The cases $\eta>1$ and $\eta\approx1$
will be considered below.

\subsection{$\eta>1$}

In this region $F(x,\eta)$ is concentrated near $x=X_{0}$, and using
(\ref{cornerlayersol}) and (\ref{P}) we have
\begin{align*}
F  &  \sim\exp\left\{  -\frac{\eta^{2}}{2\varepsilon}-\frac{\eta
}{2Dj\varepsilon^{\frac{1}{3}}}\omega^{2}\right\}  \frac{1}{\varepsilon\pi
}2^{-\frac{2}{3}}\sqrt{\frac{\eta}{Dj}}\wp(0)\\
&  =\exp\left\{  -\frac{\eta^{2}}{2\varepsilon}-\frac{\eta}{2Dj\varepsilon
}(x-X_{0})^{2}\right\}  \frac{1}{2\pi\varepsilon}\sqrt{\frac{\eta}{Dj}},\quad
x\rightarrow X_{0}%
\end{align*}
and hence, by Laplace's method,
\begin{align*}
\int\limits_{0}^{\infty}F(x,\eta)dx  &  \sim\int\limits_{-\infty}^{\infty}%
\exp\left\{  -\frac{\eta^{2}}{2\varepsilon}-\frac{\eta}{2Dj\varepsilon
}(x-X_{0})^{2}\right\}  \frac{1}{2\pi\varepsilon}\sqrt{\frac{\eta}{Dj}}dx\\
&  =\frac{1}{\sqrt{2\pi\varepsilon}}\exp\left(  -\frac{\eta^{2}}{2\varepsilon
}\right)  .
\end{align*}
This verifies (\ref{marginal}) (at least asymptotically as $\varepsilon
\rightarrow0)$ for $\eta>1.$

\subsection{$\eta\approx1$}

\bigskip For $\eta\rightarrow1$ and $x$ small we use the corner layer
expansion, i.e.,
\begin{align*}
F(x,\eta)  &  \sim\varepsilon^{-\frac{7}{6}}\frac{1}{\sqrt{2\pi}2^{\frac{1}%
{3}}D^{\frac{2}{3}}}\exp\left\{  -\frac{\eta^{2}}{2\varepsilon}+\frac
{\gamma\mu}{2D}-\frac{\gamma^{3}}{12D}\right\} \\
&  \times\frac{1}{2\pi i}\int\limits_{Br}\exp\left\{  2^{-\frac{2}{3}%
}D^{-\frac{1}{3}}\gamma\lambda\right\}  \frac{\operatorname{Ai}\left(
\lambda+2^{-\frac{1}{3}}D^{-\frac{2}{3}}\mu\right)  }{\left[
\operatorname{Ai}\left(  \lambda\right)  \right]  ^{2}}d\lambda.
\end{align*}
where $x=\mu\varepsilon^{\frac{2}{3}}$ and $\eta-1=$ $\gamma\varepsilon
^{\frac{1}{3}}.$ In the local variable $\mu,$(\ref{marginal}) becomes
\[
\int\limits_{0}^{\infty}F(x,\eta)d\mu=\varepsilon^{-\frac{7}{6}}\frac{1}%
{\sqrt{2\pi}}\exp\left\{  -\frac{\eta^{2}}{2\varepsilon}\right\}
\]
so we have to show that
\begin{equation}
\Lambda(\gamma)=2^{\frac{1}{3}}D^{\frac{2}{3}}\exp\left\{  \frac{\gamma^{3}%
}{12D}\right\}  \label{newmarginal}%
\end{equation}
where
\begin{align}
\Lambda(\gamma)  &  =\int\limits_{0}^{\infty}\frac{1}{2\pi i}\int
\limits_{Br}\exp\left\{  \left(  \frac{\mu}{2D}+2^{-\frac{2}{3}}D^{-\frac
{1}{3}}\lambda\right)  \gamma\right\}  \frac{\operatorname{Ai}\left(
\lambda+2^{-\frac{1}{3}}D^{-\frac{2}{3}}\mu\right)  }{\left[
\operatorname{Ai}\left(  \lambda\right)  \right]  ^{2}}d\lambda d\mu
\nonumber\\
&  =2^{\frac{1}{3}}D^{\frac{2}{3}}\frac{1}{2\pi i}\int\limits_{Br}%
\int\limits_{\lambda}^{\infty+i\operatorname{Im}(\lambda)}\exp\left\{
2^{-\frac{2}{3}}D^{-\frac{1}{3}}\gamma\rho\right\}  \frac{\operatorname{Ai}%
\left(  \rho\right)  }{\left[  \operatorname{Ai}\left(  \lambda\right)
\right]  ^{2}}d\rho d\lambda. \label{newmarginal1}%
\end{align}

Taking the derivative of $\Lambda$ and using \cite{AS} $\operatorname{Ai}%
^{\prime\prime}\left(  \rho\right)  =\rho\operatorname{Ai}\left(  \rho\right)
$ yields%

\begin{align}
\Lambda^{\prime}(\gamma)  &  =2^{\frac{1}{3}}D^{\frac{2}{3}}\frac{1}{2\pi
i}\int\limits_{Br}\int\limits_{\lambda}^{\infty+i\operatorname{Im}(\lambda
)}2^{-\frac{2}{3}}D^{-\frac{1}{3}}\rho\exp\left\{  2^{-\frac{2}{3}}%
D^{-\frac{1}{3}}\gamma\rho\right\}  \frac{\operatorname{Ai}\left(
\rho\right)  }{\left[  \operatorname{Ai}\left(  \lambda\right)  \right]  ^{2}%
}d\rho d\lambda\label{4.4}\\
&  =2^{-\frac{1}{3}}D^{\frac{1}{3}}\frac{1}{2\pi i}\int\limits_{Br}%
\int\limits_{\lambda}^{\infty+i\operatorname{Im}(\lambda)}\exp\left\{
2^{-\frac{2}{3}}D^{-\frac{1}{3}}\gamma\rho\right\}  \frac{\operatorname{Ai}%
^{\prime\prime}\left(  \rho\right)  }{\left[  \operatorname{Ai}\left(
\lambda\right)  \right]  ^{2}}d\rho d\lambda.\nonumber
\end{align}
Two integrations by parts give
\begin{align*}
&  \int\limits_{\lambda}^{\infty+i\operatorname{Im}(\lambda)}\exp\left\{
2^{-\frac{2}{3}}D^{-\frac{1}{3}}\gamma\rho\right\}  \operatorname{Ai}%
^{\prime\prime}\left(  \rho\right)  d\rho\\
&  =\left[  \operatorname{Ai}\left(  \lambda\right)  \right]  ^{2}\frac
{d}{d\lambda}\left[  \exp\left\{  2^{-\frac{2}{3}}D^{-\frac{1}{3}}%
\gamma\lambda\right\}  \frac{1}{\operatorname{Ai}\left(  \lambda\right)
}\right] \\
&  +\left(  2^{-\frac{2}{3}}D^{-\frac{1}{3}}\gamma\right)  ^{2}\int
\limits_{\lambda}^{\infty+i\operatorname{Im}(\lambda)}\exp\left\{
2^{-\frac{2}{3}}D^{-\frac{1}{3}}\gamma\rho\right\}  \operatorname{Ai}\left(
\rho\right)  d\rho
\end{align*}
which when used in (\ref{4.4}) leads to the differential equation
\begin{equation}
\Lambda^{\prime}(\gamma)=\frac{1}{4D}\gamma^{2}\Lambda(\gamma). \label{ODE3}%
\end{equation}

\bigskip Solving (\ref{ODE3}) yields
\begin{equation}
\Lambda(\gamma)=\Lambda_{0}\exp\left\{  \frac{\gamma^{3}}{12D}\right\}
\label{4.5}%
\end{equation}
where $\Lambda_{0}$ is a constant. To determine $\Lambda_{0}$ we let
$\gamma\rightarrow-\infty$ in (\ref{newmarginal1}). Expanding the double
integral by a combination of the Laplace and saddle point methods leads to
\[
\Lambda(\gamma)\sim D^{\frac{2}{3}}2^{\frac{1}{3}}\exp\left\{  \frac
{\gamma^{3}}{12D}\right\}  ,\quad\gamma\rightarrow-\infty.
\]
Comparing this to (\ref{4.5}) we obtain $\Lambda_{0}=D^{\frac{2}{3}}%
2^{\frac{1}{3}},$ which verifies (\ref{newmarginal}).

\subsection{The marginal distribution $M(x)$}

\bigskip To evaluate (\ref{mm}) by Laplace's method, we find where $\Psi$ and
$\Phi$ are maximal as functions of $\eta.$ We thus examine the equations
$\Psi_{\eta}=0$ and $\Phi_{\eta}=0.$

We recall from (\ref{rayII1}) that $\Phi_{\eta}=(a-b)e^{\tau}-a.$ The equation
$\Phi_{\eta}=0$ then reads
\[
e^{\tau}=\frac{\sqrt{D}(\sigma-1)}{\sqrt{D}(\sigma-1)+D^{\frac{3}{2}}%
\sigma+D\sqrt{D\sigma^{2}+\left(  \sigma-1\right)  ^{2}}}<1,\quad\text{for all
}D>0,\ \sigma>1
\]
We conclude that there is no solution to $\Phi_{\eta}=0$ for $\tau>0$, and
hence $\Phi_{\eta}<0$ in Region II.

From (\ref{ray2}) $\Psi_{\eta}=(A-B)e^{t}-A$ and consequently
\begin{equation}
\Psi_{\eta}=0\quad\Leftrightarrow\quad t=\ln\left[  \frac{1-s}{1-(D+1)s}%
\right]  \label{deriveta}%
\end{equation}
which when used in (\ref{ray3}) yields
\begin{align}
\Psi_{\eta}  &  =0\quad\Leftrightarrow\quad x=X_{1}(\eta),\quad0\leq\eta
<\frac{1}{D+1},\nonumber\\
X_{1}(\eta)  &  =-2\eta-\frac{1}{D}\left(  2D\eta-D+2\eta-2\right)  \ln\left[
\frac{1-\eta}{1-(D+1)\eta}\right]  . \label{X1}%
\end{align}

The equation $x=X_{1}(\eta)$ defines implicitly $\eta$ as a function of $x, $
$\eta=E(x).$ We introduce the function
\begin{equation}
\Psi_{1}(x)\equiv\Psi\left[  x,E(x)\right]  \label{psi1}%
\end{equation}
and from (\ref{ray3}) we get
\[
\Psi_{1}(x)=\frac{E(x)\left[  1-E(x)\right]  }{D}+\frac{D+1}{D^{2}}\left[
1-E(x)\right]  ^{2}\ln\left[  \frac{1-(D+1)E(x)}{1-E(x)}\right]  .
\]
From the defining equation%

\begin{equation}
-2E(x)+\frac{1}{D}\left[  2(D+1)E(x)-D-2\right]  \ln\left[  \frac
{1-(D+1)E(x)}{1-E(x)}\right]  =x \label{E}%
\end{equation}
we obtain the asymptotic results
\begin{align}
E(x)  &  \sim\frac{x}{D}-\frac{1}{2}\frac{x^{2}}{D}+\frac{1}{6}\frac
{D-4}{D^{2}}x^{3},\quad x\rightarrow0\label{Easympt}\\
E(x)  &  \sim\frac{1}{D+1}-\frac{D}{(D+1)^{2}}\exp\left(  -x-\frac{2}%
{D+1}\right)  ,\quad x\rightarrow\infty.\nonumber
\end{align}

\bigskip Use of Laplace' s method to evaluate the integral in (\ref{mm}) as
$\varepsilon\rightarrow$0 yields%

\[
M(x)\sim\varepsilon^{-\frac{3}{2}}K\left[  x,E(x)\right]  \sqrt{2\pi}\frac
{1}{\sqrt{-\varepsilon^{-1}\Psi_{\eta\eta}\left[  x,E(x)\right]  }}%
\exp\left\{  \frac{1}{\varepsilon}\Psi_{1}(x)\right\}
\]
and from (\ref{ray3}) after some algebra we have
\begin{align}
M(x)  &  \sim\varepsilon^{-1}\frac{\left[  1-E(x)\right]  ^{2}}{\sqrt{\Delta}%
}\exp\left\{  \frac{1}{\varepsilon}\Psi_{1}(x)\right\} \label{M1}\\
\Delta &  =\frac{2\left[  1-(D+1)E(x)\right]  \left[  1-E(x)\right]  \left[
x+2E(x)\right]  (D+1)D}{2(D+1)E(x)-D-2}\nonumber\\
&  +D\left[  D+2E(x)-2(D+1)E(x)^{2}\right]  .\nonumber
\end{align}

We can get more explicit results if $x$ is either small or large, using
(\ref{Easympt}). We obtain
\begin{gather}
M(x)\sim\varepsilon^{-1}\frac{1}{D}\left(  1-\frac{x}{D}\right)  \exp\left\{
\frac{1}{\varepsilon}\left(  -\frac{x}{D}+\frac{x^{2}}{2D^{2}}\right)
\right\}  ,\quad x\rightarrow0,\label{Mlimit}\\
M(x)\sim\varepsilon^{-1}\left[  \frac{D}{(1+D)^{2}}+\frac{2D+1}{D(1+D)^{2}%
}e^{-x-\frac{2}{D+1}}\right] \nonumber\\
\times\exp\left[  -\frac{1}{\varepsilon}\left(  \frac{x}{1+D}+\frac
{1}{(1+D)^{2}}\right)  \right]  ,\quad x\rightarrow\infty.\nonumber
\end{gather}
The first result in (\ref{Mlimit}) shows that $M(x)$ is concentrated in the
range $x=O(\varepsilon)$ and the second result is consistent with the spectral
solution to (\ref{scaled}) obtained in \cite{KT}.

\section{Acknowledgment}

The work of C. Knessl was partially supported by NSF grant DMS 99-71656. The
work of D. Dominici was supported by NSF grant DMS 99-73231, provided by
Professor Floyd Hanson. D. Dominici wish to thank him for his generous sponsorship.

\end{document}